\documentclass{amsart}

\usepackage{amsmath}
\usepackage{amssymb}
\usepackage{amsthm}
\usepackage{color}
\usepackage{a4wide}
\usepackage{hyperref}
\usepackage{cite}

\newtheorem{definition}{Definition}[section]
\newtheorem{theorem}[definition]{Theorem}
\newtheorem{lemma}[definition]{Lemma}
\newtheorem{corollary}[definition]{Corollary}
\newtheorem{proposition}[definition]{Proposition}

\newtheorem{remark}[definition]{Remark}

\newcommand{\R}{\mathbb{R}}
\newcommand{\Rn}{\mathbb{R}^n}
\newcommand{\ve}{\vert}
\newcommand{\lv}{\left\vert}
\newcommand{\rv}{\right\vert}
\newcommand{\Ve}{\Vert}
\newcommand{\lV}{\left\Vert}
\newcommand{\rV}{\right\Vert}
\newcommand{\lb}{\left\lbrace}
\newcommand{\rb}{\right\rbrace}

\newcommand{\dist}{\mathrm{dist}}
\newcommand{\diam}{\mathrm{diam}}

\def\1{{\bf 1}}

\def\nn{\nonumber}

\def\sH {{\mathcal H}}
\def\sP {{\mathcal P}}

 \def\bN {{\mathbb N}} 
  \def\bR {{\mathbb R}}

\def\R {{\mathbb R}}

\numberwithin{equation}{section}

\def\qed{{\hfill $\Box$ \bigskip}}

\def\DD{{\mathcal D}}
\def\FF{{\mathcal F}}

\def\R{{\mathbb R}}

\def\E{{\mathbb E}}

\def\P{{\mathbb P}}
\def\N{{\mathbb N}}
\def\eps{\varepsilon}

\def\wt{\widetilde}
\def\pf{\noindent{\bf Proof.} }

\numberwithin{equation}{section}
\makeatletter
\newcommand*{\rom}[1]{\expandafter\@slowromancap\romannumeral #1@}
\makeatother

\begin{document}

\title{Boundary regularity for nonlocal operators with kernels of variable orders}
\author{Minhyun Kim}
\address{Department of Mathematical Sciences, Seoul, Korea}
\email{201421187@snu.ac.kr}
\author{Panki Kim}
\address{Department of Mathematical Sciences, Seoul, Korea}
\email{pkim@snu.ac.kr}
\author{Jaehun Lee}
\address{Department of Mathematical Sciences, Seoul, Korea}
\email{hun618@snu.ac.kr}
\author{Ki-Ahm Lee}
\address{Department of Mathematical Sciences, Seoul, Korea}
\email{kiahm@snu.ac.kr}


\maketitle

\begin{abstract}
We study the boundary regularity of solutions of the Dirichlet problem for the nonlocal operator with a kernel of variable orders. Since the order of differentiability of the kernel is not represented by a single number, we consider the generalized H\"older space. We prove that there exists a unique viscosity solution of $Lu = f$ in $D$, $u=0$ in $\Rn \setminus D$, where $D$ is a bounded $C^{1,1}$ open set, and that the solution $u$ satisfies $u \in C^V(D)$ and $u/V(d_D) \in C^\alpha (D)$ with the uniform estimates, where $V$ is the renewal function and $d_D(x) = \dist(x, \partial D)$.
\end{abstract}

\tableofcontents

\section{Introduction}

In this paper, we will consider the viscosity solutions for the following Dirichlet (exterior) problem
\begin{equation}\label{e:pde}
\begin{cases}
-\phi(-\Delta)u=f &\text{in} ~ D, \\
u=0 &\text{in} ~ \Rn \backslash D,
\end{cases}
\end{equation}
where $\phi$ is in the class of functions called \textit{Bernstein function}, which contains  $\phi(\lambda)=\lambda^{\alpha}$ with $0<\alpha<1$, and $D$ is a bounded $C^{1,1}$ open set in $\Rn$. For example, if $\phi(\lambda)=\lambda^\alpha$, then $-\phi(-\Delta)=-(-\Delta)^\alpha$ is a fractional Laplacian.

We will focus on the boundary behavior of the viscosity solutions of the Dirichlet problem \eqref{e:pde} under assumptions \eqref{e:phi-wsc} and \eqref{e:exp-j} below.

\subsection{Probabilistic point of view} \label{s:prob}

The operator $-\phi (-\Delta)$ can be understood as the infinitesimal generator of subordinate Brownian motions, thus we can use probabilistic tools to study the behavior of solutions of \eqref{e:pde}.

Let $S=(S_t)_{t \ge 0}$ be a subordinator, that is, an increasing L\'evy process in $\bR$. It is known that its Laplace exponent is given by
$$ \E[e^{-\lambda S_t}] = \exp(-t\phi(\lambda)), \quad \lambda > 0, $$
where the function $\phi:(0,\infty) \rightarrow (0,\infty)$ satisfies $\displaystyle \lim_{\lambda \downarrow 0} \phi(\lambda)=0$ and
\begin{equation}\label{d:phi}
\phi(\lambda)= b\lambda + \int_{(0,\infty)} (1-e^{-\lambda x}) \mu(dx)
\end{equation}
with a drift $b \ge 0$ and a measure $\mu$ on $(0,\infty)$ satisfying $\int_{(0,\infty)} (1 \land x) \mu(dx) < \infty.$ 
It is known that the function $\phi$ of the form \eqref{d:phi} is a \textit{Bernstein function}, it means, $\phi:(0,\infty) \rightarrow (0,\infty)$ is a $C^\infty$-function satisfying
$$ (-1)^{n+1} \phi^{(n)}(\lambda) \ge 0 \quad \mbox{for all} \quad n \in \bN. $$
Here $\phi^{(n)}$ is the $n$-th derivative of $\phi$. Also, it is known that every Bernstein function can be uniquely represented by \eqref{d:phi}.

Subordinate Brownian motion $Y=(Y_t)_{t\geq 0} = (B_{S_t})_{t \geq 0}$ in $\Rn$ is a L\'evy process obtained by replacing the time of Brownian motion in $\Rn$ by an independent subordinator. Then, the characteristic exponent of $Y$ is given by $z \mapsto \phi(\ve z \ve^2)$. Also, the L\'evy measure of the process has a density $y \mapsto j(|y|)$ where $j:(0,\infty) \rightarrow (0, \infty)$ is the function given by
\begin{equation}\label{d:j}
j(r) =j_n(r)= \int_0^\infty (4\pi t)^{-n/2} e^{-\frac{r^2}{4t}} \mu(dt),
\end{equation}
and we have
\begin{equation}\label{d:Phi}
\phi(\ve z\ve^2)= \int_{\Rn \backslash \{0\}} (1- \cos ( z \cdot y)) j(|y|)dy.
\end{equation}
Let $A$ be the infinitesimal generator of $Y$. Then, by \cite[Section 4.1]{Sko} we have
\begin{equation}\label{e:SBM}
Au(x) =-\phi(-\Delta)u(x)= \int_{\Rn \setminus \lb 0 \rb} \left( u(x+y) - u(x) - {\bf 1}_{\lb \ve y \ve \le 1 \rb} y \cdot \nabla u(x) \right) j(|y|)dy.
\end{equation}
for any $u \in C^2(\R^n)$. See 
Section \ref{s:N} for the definition of function spaces and 
Section \ref{s:preliminaries} for the definition of infinitesimal generator.

Note that when $\phi(\lambda)=\lambda^\alpha$ with $0<\alpha<1$, the corresponding subordinate Brownian motion in $\Rn$ is a rotationally symmetric $2\alpha$-stable process. We also have $j(|y|)=c(n,\alpha)|y|^{-n-2\alpha}$. Thus the corresponding infinitesimal generator is the fractional Laplacian $-(-\Delta)^\alpha$.  

Now we introduce some conditions which we will impose in this paper. The first condition is \textit{weak scaling condition at the infinity} for $\phi$, that is, there exist constants $0 < \alpha_1 \le \alpha_2 < 1$ and $b_1 \geq 1$ such that
\begin{equation}\label{e:phi-wsc}
b_1^{-1} \left( \frac{R}{r} \right)^{\alpha_1} \le \frac{\phi(R)}{\phi(r)} \le b_1 \left( \frac{R}{r} \right)^{\alpha_2} \quad \mbox{for all} \quad 1 \le r \le R < \infty.
\end{equation}
The constant $1$ in above condition can be changed into other positive constant without loss of generality. Note that \eqref{d:phi} and \eqref{e:phi-wsc} imply that $b =0$ and that $\mu$ is an infinite measure. The second one is that the L\'evy density of process satisfies
\begin{equation}\label{e:exp-j}
j(r+1) \le b_2 j(r) \quad \mbox{for all} \quad r \ge 1
\end{equation}
for some constant $b_2 > 0$. \eqref{e:exp-j} is valid for any complete Bernstein function satisfying \eqref{e:phi-wsc}. See \cite[Definition 6.1]{SSV} and \cite[Theorem 13.3.5]{KSV} for details. Moreover, we also have \eqref{e:exp-j} when \eqref{e:phi-wsc} holds for any $0< r \le R < \infty$ (See \cite[Corollary 22]{BGR1}).

We will see that the \textit{renewal function} $V$ with respect to one dimensional L\'evy process is related to the boundary behavior of solutions. This function plays an important role throughout this paper. For the definition of the \textit{renewal function}, see Section \ref{s:renewal function}.

\subsection{Analytic point of view}

In analytic point of view, nonlocal operators can be defined via the Fourier transformation. For example, the fractional Laplacian is defined by
\begin{align*}
-(-\Delta)^{\sigma/2} f(x) 
&:= -(\ve \xi \ve^\sigma \hat{f})^\vee (x) \\
&= P.V. \int_{\Rn} \frac{f(y) - f(x)}{\ve y - x \ve^{n+\sigma}} \, dy \\
&= \int_{\Rn} \frac{f(y) - f(x) - \nabla f(x) \cdot (y-x) {\bf 1}_{\lb \ve y - x \ve < k\rb}}{\ve y - x \ve^{n+\sigma}} \, dy
\end{align*}
for $f \in C^\infty_c(\Rn)$ and it is well-known that
\begin{align*}
\lim_{\sigma \rightarrow 2} (2-\sigma) c(n, \sigma) (-\Delta)^{\sigma/2} f(x) = - \Delta f(x).
\end{align*}
Moreover, Caffarelli and Silvestre \cite{CS2} provided Harnack inequality and interior $C^{1,\alpha}$ regularity for fully nonlinear integro-differential equations associated with kernels comparable to that of fractional Laplacian, which remain uniform as $\sigma \rightarrow 2$. These results were generalized in \cite{KL4} and \cite{KKL} to more general integro-differential equations. These results make the theory of integro-differential operators and elliptic differential operators become unified.

The fractional Laplacian $(-\Delta)^{\sigma/2} f$ can be also thought as the normal derivative of some extension of $f$ (the Dirichlet to Neumann operator of $f$). Consider the extension problem
\begin{align*}
\begin{cases}
- \nabla (y^{1-\sigma} \nabla u) = 0 &\text{in} ~ \Rn \times (0, \infty), \\
u(x,0) = f(x) &\text{for} ~ x \in \Rn.
\end{cases}
\end{align*}
It is known in \cite{CS1} that the following holds:
\begin{align*}
(-\Delta)^{\sigma/2} f(x) = \partial_\nu u(x,0) = - \lim_{y \rightarrow 0} y^{1-\sigma} u_y(x,y),
\end{align*}
where $\partial_\nu u$ is the outward normal derivative of $u$ on the boundary $\lb y = 0\rb$.

We are interested in the operator of the form
\begin{equation}\label{e:operator}
Lu(x) = P.V. \int_{\Rn \setminus \lb 0 \rb} \left( u(x+y) - u(x) \right) j(\ve y \ve) \, dy
\end{equation}
where $j: (0, \infty) \rightarrow (0, \infty)$ is an non-increasing function satisfying 
\eqref{d:Phi}, \eqref{e:phi-wsc} and \eqref{e:exp-j}, or satisfying \eqref{e:Phi-wsc} and \eqref{e:J} in Section \ref{s:Levy process}. Let us call the function $j(|y|)$ be the \textit{kernel} of operator $L$. Note that
$Lu(x)$ is well-defined if $u \in C^2(x) \cap B(\Rn)$, where $C^2(x)$ denotes the family of all functions which are $C^2$ in some neighborhood of $x$ and $B(\Rn)$ denotes the family of all bounded functions defined on $\Rn$, and this is why we needed the assumption $0 < \alpha_1 \leq \alpha_2 < 1$. Due to the symmetry of the kernel $j(|y|)dy$, the operator can be rewritten without the principal value as
\begin{align} \label{e:operator1}
\begin{split}
Lu(x) &= \int_{\Rn \setminus \lb 0 \rb} \left( u(x+y) - u(x) - {\bf 1}_{\lb \ve y \ve \le 1 \rb} y \cdot \nabla u(x) \right) j(|y|)dy \\ 
&= \frac{1}{2}\int_{\Rn \setminus \lb 0 \rb} \left( u(x+y) + u(x-y) -2u(x) \right) j(\ve y \ve) \, dy
\end{split}
\end{align}
when $u \in C^2(x) \cap B(\Rn)$. The important point to note here is that $Lu=Au$ for $u \in C^2(\R^n)$ when $j(|y|)$ in \eqref{e:SBM} and \eqref{e:operator} are the same. In Section \ref{s:MH} we discuss the connection between two operators in \eqref{e:SBM} and \eqref{e:operator}.

We will consider the {\it viscosity solution} of $Lu = f$ in $D$. A function $u : \Rn \rightarrow \R$ which is upper (resp. lower) semicontinuous on $\overline{D}$ is said to be a {\it viscosity subsolution} (resp. {\it viscosity supersolution}) to $Lu = f$, and we write $Lu \geq f$ (resp. $Lu \leq f$) in {\it viscosity sense}, if for any $x \in D$ and a test function $v \in C^2(x)$ satisfying
$v(x) = u(x)$ and $$v(y) > u(y) \quad (\mbox{resp. }< \,), \quad y \in \Rn \setminus \lb x \rb,$$
it holds that
\begin{align*}
Lv(x) \geq f(x) \quad  (\mbox{resp.} \le).
\end{align*}
A function $u$ is said to be a {\it viscosity solution} if $u$ is both sub and supersolution.

We are going to prove the H\"older regularity of viscosity solutions of nonlocal Dirichlet problem
\begin{align} \label{e:pde2}
\begin{cases}
Lu = f &\text{in } D, \\
u = 0 &\text{in } \Rn \setminus D,
\end{cases}
\end{align}
up to the boundary using the gradient heat kernel estimates and prove higher boundary regularity using PDE tools: barriers, comparison principle, and Harnack inequality. It is important that the boundary condition in \eqref{e:pde2} is given not only on $\partial D$ but on the whole complement of $D$ because of the nonlocal character of the operator $L$. See Section \ref{s:MH} for details.

The PDE approach can be applied to nonlinear integro-differential equations. There are many literatures dealing with regularity results with PDE approach. See \cite{CS2, KL2, KL4, Bae, BK3, ROS2} and \cite{KKL}. We expect that similar results such as Harnack inequality and H\"older regularity hold for nonlinear equations with our $L$.

\subsection{History}

Over the last few decades there have been a lot of studies for the nonlocal operators, and regularity theory for nonlocal operators is one of the main areas as the one for local operators. In \cite{BL} Bass and Levin proved H\"older regularity of harmonic functions with respect to a class of pure jump Markov processes in $\Rn$, whose kernels are comparable to those of symmetric stable processes. Bass and Kassmann generalized this result to kernels with variable order in \cite{BK1, BK2}. Bass also established in \cite{Bas} the Schauder estimates for stable-like operators in $\Rn$. All these works were done by probabilistic methods.

On the other hand, in \cite{Sil} Silvestre provided a purely analytic proof of H\"older estimates for solutions to integro-differential equation. His assumptions include the case of an operator with variable orders. In \cite{CS2} Caffarelli and Silvestre generalized this result to fully nonlinear integro-differential equations associated with symmetric kernels comparable to fractional Laplacian by PDE methods. Kim and Lee, in \cite{KL2} and \cite{KL4}, extended this result to fully nonlinear integro-differential equations associated with nonsymmetric kernels. A singular regularity theory for parabolic nonlocal nonlinear equations was also established at \cite{KL3}. In \cite{Bae}, Bae proved H\"older regularity for solutions of fully nonlinear integro-differential equations with kernels of variable orders in \cite{Bae}. Bae and Kassmann in \cite{BK3} established Schauder estimates for integro-differential equation with kernels of variable orders. In \cite{KKL}, they extended the regularity results for the integro-differential operators of the fractional Laplacian type by Caffarelli and Silvestre \cite{CS2} to those for the integro-differential operators associated with symmetric, regularly varying kernels at zero.

There are relatively fewer results concerning boundary regularity of solutions of Dirichlet problem. For the boundary regularity for local operators, see \cite{DL}. Kim and Lee proved regularity up to the boundary for the fractional heat flow in \cite{KL1}. The boundary regularity up to the boundary is well-known for the fractional Laplacian, and for fully nonlinear integro-differential equations, when $D$ is a bounded $C^{1,1}$ domain. See \cite{ROS1, ROS2}. Ros-Oton and Serra also proved the similar result when $D$ is a bounded $C^{1, \alpha}$ or $C^1$ domain in \cite{ROS3}. However, there is no boundary regularity result for the operators with kernels having variable orders.

\subsection{Notation} \label{s:N}

In this paper, we denote $a \land b = \min \{a,b\}$ and $a \lor b = \max\{a,b\}$. For any nonnegative functions $f$ and $g$, $f(r) \asymp g(r)$ for $r > 0$ (resp. $0<r \le r_0$) means that there is a constant $ c \ge 1$ such that $c^{-1} f(r) \le g(r) \le c f(r)$ for $r>0$ (resp. $0<r \le r_0$). We call $c$ the \textit{comparison constant} of $f$ and $g$.  We also denote $B(x,r):= \{ y \in \R^n : |x-y|<r  \}$ for the open ball and $d_D(x) := \dist (x, D^c)$ for the distance between $x \in D$ and $D^c$. For $n \ge 1$, let $\omega_n= \int_{\R^n} \1_{ \{|y| \le 1\} } dy$ be the volume of $n$-dimensional ball.

We denote by $C(D)$ the Banach space of bounded and continuous functions on $D$, equipped with the supremum norm $\Ve f \Ve_{C(D)} := \sup_{x \in D} \ve f(x) \ve$, and denote by $C^k(D), k \geq 1$, the Banach space of $k$-times continuously differentiable functions on $D$, equipped with the norm $\Ve f \Ve_{C^k(D)} := \sum_{\ve \gamma \ve \leq k} \sup_{x\in D} \ve D^\gamma f(x) \ve$. Also, denote $C_0(D):= \{ u \in C(D): u \mbox{ vanishes at the boundary of } D \}$. For $x \in \R^n$, define $C^{1}(x)$ as the collection of functions which are $C^{1}$ in some neighborhood of $x$. Similarly, we define $C^2(x)$, $C^{1,1}(x)$, etc.  For $0<\alpha<1$, the H\"older space $C^\alpha(\R^n)$ is defined as
\begin{align} \label{e:Holder}
C^\alpha(\R^n) := \lb f \in C(\R^n) ~ | ~ \Ve f \Ve_{C^\alpha(\Rn)} < \infty \rb,
\end{align}
equipped with the $C^\alpha$-norm 
$$ \Ve f \Ve_{C^\alpha(\R^n)} := \Ve f \Ve_{C(\R^n)} + \sup_{x,y \in \Rn, x\neq y} \frac{|f(x)-f(y)|}{|x-y|^\alpha}.$$
Also, for given open set $D \subset \Rn$ we define $C^\alpha(D)$ by
$$ C^\alpha(D) := \lb f \in C(D) ~| ~ \Ve f \Ve_{C^{\alpha}(D)} < \infty \rb $$
with the norm
$$ \Ve f \Ve_{C^\alpha(D)} := \Ve f \Ve_{C(D)} + \sup_{x,y \in D, x\neq y} \frac{|f(x)-f(y)|}{|x-y|^\alpha}.$$
For given function $h : (0, \infty) \rightarrow (0, \infty)$, we define Generalized H\"older space $C^h(D)$ for bounded open set $D$ by
\begin{align} \label{e:gen Holder}
C^h(D) := \left\{f \in C(D) ~ | ~ \Ve f \Ve_{C^h(D)} < \infty \rb,
\end{align}
equipped with the norm
$$\Ve f \Ve_{C^h(D)} := \Ve f \Ve_{C(D)} + \sup_{x,y \in D, x\neq y}\frac{|f(x)-f(y)|}{h(|x-y|)}.
$$
We define seminorm $[\, \cdot \,]_{C^h(D)}$ by 
$$ [f]_{C^h(D)} := \sup_{x,y \in D, x\neq y}\frac{|f(x)-f(y)|}{h(|x-y|)}. $$
We denote the diameter  of $D$ by  diam$(D)$.
Note that if $h_1 \asymp h_2$ in $0<r \le $ diam$(D)$, $\Ve \cdot \Ve_{C^{h_1}(D)}$ and $\Ve \cdot \Ve_{C^{h_2}(D)}$ are equivalent and $C^{h_1}(D) = C^{h_2}(D)$.

We say that $D\subset \Rn$ (when $n\ge 2$) is a $C^{1,1}$ open set if there exist a localization radius $ R_0>0 $ and a constant $\Lambda>0$ such that for every $z\in\partial D$ there exist a $C^{1,1}$-function $\varphi=\varphi_z: \R^{n-1}\to \R$ satisfying $\varphi(0)=0$, $\nabla\varphi (0)=(0, \dots, 0)$, $\| \nabla\varphi \|_\infty \leq \Lambda$, $| \nabla \varphi(x)-\nabla \varphi(w)| \leq \Lambda |x-w|$ and an orthonormal coordinate system $CS_z$ of  $z=(z_1, \cdots, z_{n-1}, z_n):=(\wt z, \, z_n)$ with origin at $z$ such that $ D\cap B(z, R_0 )= \{y=({\tilde y}, y_n) \in B(0, R_0) \mbox{ in } CS_z: y_n > \varphi (\wt y) \}$. 
The pair $( R_0, \Lambda)$ will be called the $C^{1,1}$ characteristics of the open set $D$.
Note that a $C^{1,1}$ open set $D$ with characteristics $(R_0, \Lambda)$ can be unbounded and disconnected, and the distance between two distinct components of $D$ is at least $R_0$.
By a $C^{1,1}$ open set  in $\R$ with a characteristic $R_0>0$, we mean an open set that can be written as the union of disjoint intervals so that the {infimum} of the lengths of all these intervals is {at least $R_0$} and the {infimum} of the distances between these intervals is  {at least $R_0$}.

\subsection{Main theorems}

The main results of this paper are the existence and the uniqueness of the viscosity solution $u$ of \eqref{e:pde}, the generalized H\"older regularity estimates of such solution $u$ and the regularity of the quotient $u\phi(d_D^{-2})$ up to the boundary.

The boundary estimate for nonlinear PDE has been studied for a long time, where the solution behaves as a linear function. See \cite{CC} and references therein. For the degenerate or singular PDE, \cite{KL3}, it has been proved that the solution behaves in various ways just as that of the fractional Laplace equation. In \cite{ROS1}, Ros-Oton and Serra applied the known techniques for local operators to fractional Laplacian, which has a nice scaling invariance and a simple barrier of the form $x_n^{\alpha}$.  On the other hand,  our $\phi$ has only a weak scaling condition at infinity and it has a general form which allows nontrivial boundary behavior different from $x_n^{\alpha}$.
In this paper, we  track down $u$ in every scale to find scaling invariant uniform estimates only with the weak scaling condition at infinity. We also construct the renewal function, $V(\cdot)$, of the ladder height process defined at \eqref{e:V} to overcome the lack of a simple barrier.
In addition, we provide the existence and uniqueness theory for given Dirichlet problem by utilizing the concept of viscosity solution.

The first result is the H\"older estimates up to the boundary of solutions of the Dirichlet problem \eqref{e:pde}. Unlike the case of the fractional Laplacian, it is inappropriate to represent H\"older regularity as a single number since kernel in \eqref{e:operator} has variable orders. Therefore it is natural to consider a generalized H\"older space.

\begin{theorem}[H\"older estimates up to the boundary] \label{t:est-u}
Assume that $D$ is a bounded $C^{1,1}$ open set in $\Rn$, and $\phi$ is a Bernstein function satisfying \eqref{e:phi-wsc} and \eqref{e:exp-j}. If $f \in C (D)$, then there exists a unique viscosity solution $u$ of \eqref{e:pde2} and $u \in C^{\overline{\phi}}(D)$. Moreover, we have
$$ \Ve u \Ve_{C^{\overline{\phi}} (D)} \leq C \Ve f \Ve_{C(D)}, $$
where $\overline{\phi}(r) := \phi(r^{-2})^{-1/2}$, for some constant $C > 0$ depending only on $n, D$, and $\phi$.
\end{theorem}

We will prove Theorem \ref{t:est-u} using the potential operator, which is the inverse of the operator $L$, and the estimates on the transition density and its spatial derivatives, see Section \ref{s:H} for details. In whole space $\R^n$, estimates on any order of spatial derivatives of the transition density are known. Based on these estimates, Bae and Kassmann established Schauder estimates for the integro-differential operators with kernels of variable orders in \cite{BK3}. However, in a bounded $C^{1,1}$ open set, estimates on the first order derivative of the transition density are only known. Higher order regularities up to the boundary require further research in future.

It is well known that $\bar{\phi}$ is comparable to renewal function $V$ (see Section \ref{s:renewal function}.) Thus any solution $u$ of Dirichlet problem \eqref{e:pde} is in $C^V$ up to the boundary by Theorem \ref{t:est-u}. Hence it is of importance to study the regularity of $u/V(d_D)$ up to the boundary. The following is our second main result.

\begin{theorem}[Boundary estimates] \label{t:est-u/V}
Assume that $D$ is a bounded $C^{1,1}$ open set in $\Rn$, and $\phi$ is a Bernstein function satisfying \eqref{e:phi-wsc} and \eqref{e:exp-j}. If $f \in C (D)$ and $u$ is the viscosity solution of \eqref{e:pde2}, then $u / V(d_D) \in C^\alpha(D)$ and 
\begin{align*}
\lV \frac{u}{V(d_D)} \rV_{C^\alpha(D)} \leq C \Ve f \Ve_{C(D)}
\end{align*}
for some constants $\alpha > 0$ and $C > 0$ depending only on $n,D$, and $\phi$.
\end{theorem}

One of the methods proving the above result follows the standard argument of Krylov in \cite{Kry}. In the other words, we are going to control the oscillation of the function $u \phi(d_D^{-2})^{1/2}$ near the boundary using barriers, comparison principle, and the Harnack inequality. However, the construction of barriers are highly nontrivial. The difficulty mainly comes from the fact that the operator \eqref{e:operator} is not scale-invariant.

In fact, we will prove Theorems \ref{t:est-u} and \ref{t:est-u/V} for a little more general operators including $-\phi(-\Delta)$. In section \ref{s:preliminaries} we will state the generalization of these theorems, and we collect some known results about the renewal function $V$. We will prove Theorem \ref{t:est-u} in Section \ref{s:H}, and Theorem \ref{t:est-u/V} in section \ref{s:bdry reg}.

\section{Preliminaries} \label{s:preliminaries}

The operators we consider in this paper coincides with infinitesimal generators of isotropic unimodal L\'evy processes for $C^2(\R^n)$ functions. Thus, in Section \ref{s:Levy process} we first explain the definitions and properties of L\'evy processes, and some related concepts. Then we introduce some additional conditions that will be needed in this paper. With these concepts, we state Theorems \ref{t:1} and \ref{t:2}, which are generalized version of Theorems \ref{t:est-u} and \ref{t:est-u/V}. Throughout this paper, we prove Theorems \ref{t:1} and \ref{t:2}.

Next, in Section \ref{s:renewal function} we will define the renewal function $V$, which will be act as a barrier, and record some properties of renewal function.

\subsection{L\'evy processes} \label{s:Levy process}

Let $X=(X_t,\P^x, t \ge 0, x \in \R^n)$ be a L\'evy process  in $\Rn$ defined on the probability space $(\Omega, \FF, \P^x)$ with $\P^x(X_0=x)=1$. For the precise definition of L\'evy process, see \cite[Definition 1.5]{Sat}. Note that  $\P^x(X_t \in A)=\P^0(X_t+x \in A)$. By L\'evy-Khintchine formula, the characteristic exponent of L\'evy process is given by
$$ \E^0[e^{i z \cdot X_t}] = e^{t\Phi(z)}, \quad z \in \Rn, $$
where
$$ \Phi(z) = -\frac{1}{2}z \cdot Uz + i\gamma \cdot z + \int_{\Rn} \big( e^{iz \cdot x} -1 -iz \cdot x {\bf 1}_{\{|x| \le 1 \} } \big) J(dx) $$
with an $n \times n$ symmetric nonnegative-definite matrix $U = (U_{ij})$, $\gamma \in \Rn$ and a measure $J(dx)$ on $\Rn \backslash \{0\}$ satisfying
$$ \int_{\Rn \backslash \{0\}} \big( 1 \land |x|^2 \big) J(dx) < \infty. $$
%The triplet $(A,J,\gamma)$ determines the L\'evy process $X$ and it is called \textit{generating triplet} of $X$.
Let $(P_t)_{t \ge 0 }$ be a transition semigroup for $X$, it means that 
$$P_t f(x) := \E^x[f(X_t)]= \E^0[f(x+X_t)]. $$
Now, define \textit{the infinitesimal generator} $A$ of $X$ by
$$Au(x) := \lim_{t \downarrow 0} \frac{P_t u(x) - u(x)}{t}$$
if the limit exists.
%We denote $\DD(A)$ be the domain of the operator $A$.%
By \cite[Section 4.1]{Sko}, $Au$ is well-defined for $u \in C^2(\Rn)$ and represented by
$$Au(x) = \frac{1}{2} \sum_{i,j=1}^n U_{ij}\partial_{ij}u(x) + \sum_{i=1}^n \gamma_i\partial_i u(x) + \int_{\Rn \setminus \lb 0 \rb} \left( u(x+y) - u(x) - {\bf 1}_{\lb \ve y \ve \le 1 \rb} y \cdot \nabla u(x) \right) J(dy). $$

Throughout this paper, we will assume that $X$ is an isotropic unimodal pure jump L\'evy process with an infinite L\'evy measure, that is, $U=0$, $\gamma=0$ and $J(dy)$ is an infinite measure with an isotropic density $J(|y|)dy$, where $r \mapsto J(r)$ is non-increasing. Under these assumptions, $X$ possesses transition density $p:(0,\infty) \times \R_+ \rightarrow \R_+$ satisfying
$$P_t f(x) = \E^x[f(X_t)] = \int_{\Rn} f(y) p(t,|x-y|)dy $$
and characteristic exponent $\Phi: \R^n \rightarrow \R_+$ is an isotropic function. From now on, we regard isotropic functions $J$ and $\Phi$ as functions on $\R_+$.

For every open subset $D \subset \R^n$, let $\tau_D:= \inf \{ t>0 : X_t \notin D \}$ be the first exit time of $D$ by $X$. We define subprocess $X^D=(X_t^D)_{t \ge 0}$, which is called \textit{the killed process of X upon $D$}, by $X_t^D=X_t$ when $t < \tau_D$ and $X_t^D=\partial$ when $t \le \tau_D$ where $\partial$ is a cemetery point. Since $X$ has the transition density, $X^D$ also possesses the transition density $p_D(t,x,y)$ with 
$$p_D(t,x,y) = p(t,|x-y|)-\E^x[p(t-\tau_D,|X_{\tau_D}-y|);\tau_D<t],$$
and its transition semigroup $(P^D_t)_{t \ge 0}$ is represented by
$$P^D_t f(x) : = \E^x[f(X^D_t)] = \int_{D}f(y)p_D(t,x,y) \, dy. $$
% As the L\'evy processes, $Lu(x)$ coincides with \eqref{e:SBM} when $u \in C^2_0(D)$.

Now we are ready to introduce main assumptions in this paper. Note that, under settings above, the infinitesimal generator can be rewritten as
\begin{align}\label{d:A}
Au(x) &= \frac{1}{2} \int_{\Rn \setminus \lb 0 \rb} \left( u(x+y)+ u(x-y) - 2u(x) \right)J(|y|)dy
\end{align}
for $u \in C^2(\Rn)$. Moreover, it is known in \cite[Lemma 2.6]{BLM} that \eqref{d:A} still holds for $u \in C^2(x) \cap C_0(\Rn)$.
Recall that the operator $L$ in \eqref{e:operator} with  kernel $J(|y|)$ is represented as
\begin{align}
\begin{split}\label{d:L}
Lu(x) &= \frac{1}{2} \int_{\Rn \setminus \lb 0 \rb} \left( u(x+y)+ u(x-y) - 2u(x) \right)J(|y|)dy
\end{split}
\end{align}
for $u \in C^2(x) \cap B(\Rn)$ since $J$ is symmetric. We record that $Au(x)=Lu(x)$ for any $u \in C^{2}(x) \cap C_0(\Rn)$ for the next use. 

We first assume that the characteristic exponent $\Phi$ satisfies weak scaling condition with constants $a_1 \ge 1$ and $0<\alpha_1 \le \alpha_2<1$ so that
\begin{equation}\label{e:Phi-wsc}
a_1^{-1} \left( \frac{R}{r} \right)^{2\alpha_1} \le \frac{\Phi(R)}{\Phi(r)} \le a_1 \left( \frac{R}{r} \right)^{2\alpha_2} ~~ \text{for all} ~~ 1 < r \le R \le \infty. 
\end{equation}
We also assume that 
%$X$ is isotropic unimodal pure jump L\'evy process with L\'evy measure $J(|y|)dy$, and satisfies that there exists a constant $a_2 >0$ such that
The L\'evy measure of  the isotropic unimodal pure jump L\'evy process $X$ has the density $y \to J(|y|)$ and it satisfies that there exists a constant $a_2 >0$ such that
\begin{equation}\label{e:J}
J(r+1) \le a_2 J(r) \,\, \mbox{for all} \,\,r>0, \,\, \mbox{and} \quad r \mapsto  -\frac{J'(r)}{r} \quad \mbox{is non-increasing}.
\end{equation}

\noindent Let
\begin{equation*}
\varphi(r): = \frac{J(1)}{J(r) r^n}. 
\end{equation*}
By \cite{BGR1}, for any $c>0$ we have $\Phi(r^{-1})^{-1} \asymp \varphi(r)$ in $0< r \leq c$ with comparison constant depending only on $c$ and $n$. Thus, there exists a constant $a_3=a_3(n,a_1) \geq 1$ such that
\begin{equation}\label{e:varphi-wsc}
a_3^{-1} \left( \frac{R}{r} \right)^{2\alpha_1} \le \frac{\varphi(R)}{\varphi(r)} \le a_3 \left( \frac{R}{r} \right)^{2\alpha_2} ~~ \text{for all} ~~ 0 < r \le R \le 1,
\end{equation}
where $\alpha_1$ and $\alpha_2$ are constants in \eqref{e:Phi-wsc}. 
%Note that when $X$ is unimodal, \eqref{e:Phi-wsc} is equivalent to \eqref{e:varphi-wsc}. Also, 
Note that
\eqref{e:varphi-wsc} implies that $\varphi(r) \le c r^{2\alpha_1}$ for $r \le 1$, so by definition of $\varphi$ we see that $J(|y|)dy$ is an infinite measure.

We say that $D\subset \R^d$ (when $d\ge 2$) is  a $C^{1,1}$ open set with $C^{1, 1}$ characteristics $(R_0, \Lambda)$ if there exist a localization radius $ R_0>0 $ and a constant $\Lambda>0$ such that for every $z\in\partial D$ there exist a $C^{1,1}$-function $\varphi=\varphi_z: \R^{d-1}\to \R$ satisfying $\varphi(0)=0$, $\nabla\varphi (0)=(0, \dots, 0)$, $\| \nabla\varphi \|_\infty \leq \Lambda$, $| \nabla \varphi(x)-\nabla \varphi(w)| \leq \Lambda |x-w|$ and an orthonormal coordinate system $CS_z$ of  $z=(z_1, \cdots, z_{d-1}, z_d):=(\wt z, \, z_d)$ with origin at $z$ such that $ D\cap B(z, R_0 )= \{y=({\tilde y}, y_d) \in B(0, R_0) \mbox{ in } CS_z: y_d > \varphi (\wt y) \}$. 
The pair $( R_0, \Lambda)$ will be called the $C^{1,1}$ characteristics of the open set $D$.
Note that a bounded $C^{1,1}$ open set $D$ with characteristics $(R_0, \Lambda)$ can be  disconnected, and the distance between two distinct components of $D$ is at least $R_0$.
By a $C^{1,1}$ open set  in $\R$ with a characteristic $R_0>0$, we mean an open set that can be written as the union of disjoint intervals so that the {infimum} of the lengths of all these intervals is {at least $R_0$} and the {infimum} of the distances between these intervals is  {at least $R_0$}.

Now, consider the following Dirichlet (exterior) problem on a bounded $C^{1,1}$ open set $D \subset \R^n$:
\begin{equation}\label{e:pde1}
\begin{cases}
Lu=f &\text{in} ~ D , \\
u=0 &\text{in} ~ \Rn \backslash D ,
\end{cases}
\end{equation}
where $L$ is the operator in \eqref{d:L}, which coincides with \eqref{e:pde2} when the process $X$ is a subordinate Brownian motion. We will prove the following theorems, which contain Theorem \ref{t:est-u} and \ref{t:est-u/V} (See Remark \ref{r:SBM} below), in Sections \ref{s:H} and \ref{s:bdry reg}, respectively. 

\begin{theorem} [H\"older estimates up to the boundary] \label{t:1} 
Assume that $D$ is a bounded $C^{1,1}$ open set in $\Rn$, and $X$ is an isotropic pure jump L\'evy process satisfying \eqref{e:Phi-wsc} and \eqref{e:J}.  If $f \in C(D)$, then there exists a unique viscosity solution $u$ of \eqref{e:pde1} and $u \in C^{\overline{\phi}}(D)$. Moreover, we have
$$ \Ve u \Ve_{C^{\overline{\phi}} (D)} \leq C \Ve f \Ve_{C(D)}, $$
where $\overline{\phi}(r) := \varphi(r)^{1/2}$, for some constant $C > 0$ depending only on $n, D$, and $\Phi$.
\end{theorem}

\begin{theorem} [Boundary estimates] \label{t:2}
Assume that $D$ is a bounded $C^{1,1}$ open set in $\Rn$, and $X$ is an isotropic pure jump L\'evy process satisfying \eqref{e:Phi-wsc} and \eqref{e:J}. If $f \in C (D)$ and $u$ is the viscosity solution of \eqref{e:pde1}, then $u / V(d_D) \in C^\alpha(D)$ and 
\begin{align*}
\lV \frac{u}{V(d_D)} \rV_{C^\alpha(D)} \leq C \Ve f \Ve_{C(D)}
\end{align*}
for some constants $\alpha>0$ and $C>0$ depending only on $n,D$, and $\Phi$.
\end{theorem}

\noindent In the next remark, we explain that assumptions in Theorem \ref{t:est-u} and Theorem \ref{t:est-u/V} imply assumptions in Theorem \ref{t:1} and Theorem \ref{t:2}.

\begin{remark}\label{r:SBM}
	When $X$ is a subordinate Brownian motion satisfying \eqref{e:phi-wsc} and \eqref{e:exp-j}, we have \eqref{e:Phi-wsc} by using $\Phi(r)=\phi(r^2)$ and \eqref{e:phi-wsc}. We also have that by \eqref{d:j}
$$
	J(r) = J_n(r) = \int_0^\infty (4\pi t)^{-n/2} e^{-\frac{r^2}{4t}} \mu(dt).
$$
 Thus $J(r)$ is decreasing. Also, differentiating above equation we obtain
$$ -\frac{J'_n(r)}{r}= 2\pi J_{n+2}(r), \quad r>0, $$
so $-\frac{J(r)}{r}$ is decreasing. Therefore, \eqref{e:J} holds.

Note that by \cite[Corollary 23]{BGR1} and \eqref{e:Phi-wsc} we have $\varphi(r) \asymp \Phi(r^{-1})^{-1}$. Using this and $\Phi(r)=\phi(r^2)$, both $\overline{\phi}$'s in Theorem \ref{t:est-u} and Theorem \ref{t:1} are comparable. Therefore, two $C^{\overline{\phi}}(D)$ norms are equivalent.
\end{remark}

\subsection{Renewal function} \label{s:renewal function}

Let $Z = (Z_t)_{t \geq 0}$ be an one-dimensional L\'evy process with characteristic exponent $\Phi(\ve z \ve)$ and $M_t := \sup \{ Z_s : 0 \le s \le t\}$ be the supremum of $Z$. Let $L=(L_t)_{t \ge 0}$ be a local time of $M_t - Z_t$ at 0, which satisfies
$$ L_t = \int_0^t {\bf 1}_{\{M_t = Z_t\}}(s)ds. $$
Note that since $t \mapsto L_t$ is non-decreasing and continuous with probability 1, we can define the right-continuous inverse of $L$ by
$$L^{-1}(t) := \inf \{s>0 : L(s) >t \}. $$
The mapping $t \mapsto L^{-1}(t)$ is non-decreasing and right-continuous a.s. The process $L^{-1}=(L^{-1}_t)_ {  t \ge 0  }$ with $L^{-1}_t = L^{-1}(t)$ is called \textit{the ascending ladder time process} of $Z$. \textit{The ascending ladder height process} $H=(H_t)_{t \ge 0}$ is defined as 
$$  H_t := \begin{cases} M_{L^{-1}_t} (= Z_{L^{-1}_t} ) \quad &\mbox{if} \quad L_t^{-1} < \infty, \\  \infty \quad &\mbox{otherwise}. \end{cases} $$
\noindent  (See \cite{Fri} for details.) Define the renewal function of the ladder height process $H$ with respect to $\Phi$ by 
\begin{align} \label{e:V}
V(x) = \int_0^\infty \P(H_s \le x ) ds, \quad x \in \bR.
\end{align} 
It is known that $V(x) =0 $ if $ x\le 0$, $V(\infty)=\infty$ and $V$ is strictly increasing, differentiable on $[0, \infty)$. So, there exists the inverse function $V^{-1}:[0,\infty) \rightarrow [0,\infty)$. 
%Also, $V$ is subadditive, in other words,
%\begin{equation*}
%V(x+y) \le V(x) + V(y) , \quad x,y \in \bR.
%\end{equation*}

In the following lemma we collect some basic scaling properties of renewal function in \cite{BGR1} and \cite{BGR2}. 
\begin{lemma}\label{l:V0}
	For any $c>0$, There exist constants $C_i(c)= C_i(c,n,a_1,\alpha_1,\alpha_2)>0$ for $i=1,2,3$ such that
	\begin{equation}\label{e:V-asymp}
	C_1^{-1}\varphi(r) \le V(r)^2 \le  C_1\varphi(r),\quad 0<r \le c,
	\end{equation}
	\begin{equation}\label{e:V-wsc}
	C_2^{-1} \left( \frac{R}{r} \right)^{\alpha_1} \le \frac{V(R)}{V(r)} \le C_2 \left( \frac{R}{r} \right)^{\alpha_2}, \quad 0<r\le R \le c \quad \mbox{and}
	\end{equation}
	\begin{equation}\label{e:V-inv-wsc}
C_3^{-1} \left(\frac{T}{t} \right)^{1/\alpha_2} \le \frac{V^{-1}(T)}{V^{-1}(t)} \le C_3 \left(\frac{T}{t} \right)^{1/\alpha_1}, \quad 0<t \le T < V(c).
	\end{equation}
\end{lemma}
\pf	By \cite[Corollary 3]{BGR1} and \cite[Proposition 2.4]{BGR2}, we have
\begin{equation*} 
(V(r))^{-2} \asymp \Phi(r^{-1}), \quad r>0.
\end{equation*}
with comparison constant depending only on $n$. Combining this with $\Phi(r^{-1})^{-1} \asymp \varphi(r)$ in $0<r\le c$, we conclude \eqref{e:V-asymp}.

By \eqref{e:V-asymp} and \eqref{e:varphi-wsc} we have \eqref{e:V-wsc}. Using \cite[Remark 4]{BGR1}, we also obtain the weak scaling property of the inverse function in \eqref{e:V-inv-wsc}. \qed

The most important property of renewal function in this paper is the following: $w(x):=V(x_n)$ is a solution of the following Dirichlet problem :
\begin{equation}\label{e:pde-V}
\begin{cases}
Lw=0 &\text{in} \quad  \Rn_+, \\
w=0 &\text{in} \quad  \Rn \backslash \Rn_+,
\end{cases}
\end{equation}
where $L$ is of the form \eqref{d:L} and $\bR_+^n := \left\{x=(x_1,...,x_n) \in \bR^n ~|~ x_n >0 \right\}$ is upper half plane (see \cite[Theorem 3.3]{GKK}). 

The following estimates for derivatives of $V$ are in \cite[Proposition 3.1]{GKK} and \cite[Theorem 1.2]{KR1}.

\begin{lemma}\label{l:V}
Assume $X$ is an isotropic pure jump L\'evy process satisfying \eqref{e:Phi-wsc} and \eqref{e:J}. Then $r \mapsto V(r)$ is twice-differentiable for any $r>0$. Moreover, for any $c>0$ there exists a constant $C(c)=C(c,n,a_1,\alpha_1, \alpha_2)>0$ such that 
\begin{equation}
\label{e:V-diff}
|V''(r)| \le C \frac{V'(r)}{r \land c}, \quad V'(r) \le C \frac{V(r)}{r \land c}.
\end{equation}
%In this section, we will use $C_1(C)$(or $C_1(\tilde{C})$) where $C$ and $\tilde{C}$ are constants in Definition \ref{d:psi}.
\end{lemma}
\noindent We are going to utilize the space $C^V(D)$ in Section \ref{s:H} and adopt $V(d_D)$ as a barrier in Section \ref{s:bdry reg}.

\section{H\"older regularity up to the boundary} \label{s:H}

In this section, we give the proof of Theorem \ref{t:1}. First we introduce the following Dirichlet heat kernel estimates from \cite[Corollary 1.6]{CKS} and \cite[Thoerem 1.1 and 1.2]{KR2}. We reformulate here for the usage of our proofs.

\begin{theorem}\label{t:est-p}
	Let $X$ be an isotropic unimodal L\'evy process satisfying \eqref{e:Phi-wsc} and \eqref{e:J}. Let $D \subset \Rn$ be a bounded $C^{1,1}$ open set satisfying $\diam(D) \le 1$ and $p_D(t,x,y)$ be the Dirichlet heat kernel for $X$ on $D$. Then $x \mapsto p_D(t,x,y)$ is differentiable  for any $y \in D, t >0$, and there exist constants $C_i=C_i(n,D,a_1,a_2,\alpha_1,\alpha_2,\Phi(1))>0$, $i=1,\dots,4$ satisfying the following estimates:
	\begin{enumerate}
		\item[(a)] For any $(t,x,y) \in (0,1] \times D \times D$,
		\begin{align*}
		p_D(t,x,y) \le C_1 \left( 1 \land \frac{V(d_D(x))}{t^{1/2}} \right)  \left( 1 \land \frac{V(d_D(y))}{t^{1/2}} \right)  p\left( t, |x-y|/4 \right)
		\end{align*}
		and
		$$|\nabla_x p_D(t,x,y) | \le C_2 \left[ \frac{1}{d_D(x) \land 1} \lor \frac{1}{V^{-1}(\sqrt t)} \right] p_D(t,x,y). $$
		\item[(b)] For any $(t,x,y) \in [1,\infty) \times D \times D$,
		$$ p_D(t,x,y) \le C_3 e^{-\lambda_1 t} V(d_D(x)) V(d_D(y)) $$
		and 
		$$|\nabla_x p_D(t,x,y) | \le C_4 \left[ \frac{1}{d_D(x) \land 1} \lor \frac{1}{V^{-1}(1)} \right] p_D(t,x,y), $$
	\end{enumerate}
	where $-\lambda_1 = -\lambda_1(n,a_1,a_2,\alpha_1,\alpha_2,\Phi(1))<0$ is the largest eigenvalue of the generator of $X^{B(0,1)}$.
	
\end{theorem} 
In the estimates of Theorem \ref{t:est-p}, we used $d_D(x) \lor d_D(y) \le \diam(D) \le 1$, $V(r) \asymp \varphi(r)^{1/2}$ in $0<r\le 1$ and $\frac{1}{V^{-1}(\sqrt t)} \asymp \varphi^{-1}(t)$ to reformulate theorems in our references. In addition, estimates in \cite[Corollary 1.6]{CKS} 
are of the form
$$ p_D(t,x,y) \le ce^{-\lambda(D)t} V(d_D(x))V(d_D(y))$$
where $-\lambda(D)<0$ is the largest eigenvalue of the generator of $X^{D}$. Using \cite[(6.4.14) and Lemma 6.4.5]{FOT}, we have $\lambda(D)=\inf \{ \int_{\R^n} -Lu(x)u(x)dx \, | \, \Vert u \Vert_2 =1, \mbox{supp}(u) \subset D   \}$, thus we can obtain $\lambda_1 \le \lambda(D)$. This implies heat kernel estimates in Theorem \ref{t:est-p}(b). 

Without loss of generality, we will always assume $\diam(D) \le 1$ in this paper.

\subsection{Potential operator for the killed process of subordinate Brownian motion}

In this subsection, we assume that $D \subset \R^n$ is a bounded $C^{1,1}$ open set with diam$(D)\le 1$ and $X$ is a L\'evy process satisfying \eqref{e:Phi-wsc} and \eqref{e:J}, which are conditions in Theorem \ref{t:est-p}. We define the \textit{Green function of $X^D$} by
$$ G^D (x,y) = \int_0^\infty p_D(t,x,y) dt $$
for $x,y \in D$ with $x \neq y$. Note that by Theorem \ref{t:est-p}(b), $G^D(x,y)$ is finite for any $x \neq y$.

We define a potential operator $R^D$ for $X^D$ as
\begin{equation}\label{d:R}
R^D f(x) := \int_0^\infty \int_{D} p_D(t,x,y) f(y) dy dt. 
\end{equation}
Using definitions of $P_t^D$ and $G^D$, we also have
\begin{align} \label{e:R^D}
R^Df(x) = \int_{D 
	\backslash \{x\}} G^D(x,y) f(y) dy = \int_0^\infty P_t^D f(x) dt.
\end{align}
In the next subsection, we will see that $R^D$ acts as the inverse of $-A$. 

First we will prove interior H\"older estimate of $R^D f$. For the next usage, we prove the following proposition for the functions in $L^\infty(D)$.

\begin{proposition}\label{p:R}
  	For any $f \in L^\infty(D)$ and any ball $B(x_0,r) \subset D$ satisfying $d_D(x_0) \le 2r$, we have $R^Df \in C^V(B/2)$ and there is a constant $C=C(n,a_1,a_2,\alpha_1,\alpha_2,D,\Phi(1))>0$ satisfying 
	\begin{equation}
	\label{e:R}
	\Vert R^Df \Vert_{C^V(B/2)} \le C \left( \Vert f \Vert_{L^\infty(D)} + \Vert R^Df \Vert_{C(B)} \right)
	\end{equation}
Here we have denoted $B=B(x_0,r)$ and $B/2=B(x_0,r/2)$.
\end{proposition}

\pf We have $\ve x-y \ve < r$ for any $x, y \in B/2$. Thus, we have
\begin{align*}
[R^Df]_{C^V(B/2)} 
&\leq \sup_{\ve h \ve \le r} \sup_{x \in B/2} \frac{\ve R^Df(x+h) - R^Df(x) \ve}{V(\ve h \ve)} \\
&\leq \sup_{\ve h \ve \le r} \int_0^\infty \sup_{x\in B/2} \frac{\ve P_s^Df(x+h) - P_s^Df(x) \ve}{V(\ve h \ve)} \, ds \\
&\leq \sup_{|h| \le r} \left( \int_0^{V(|h|)V(r)} + \int_{V(|h|)V(r)}^{V(r)^2} + \int_{V(r)^2}^\infty \right) \sup_{x\in B/2} \frac{|P_s^D f(x+h)-P_s^D f(x)|}{V(|h|)} \, ds \\
&=: \sup_{|h| \le r} \left( {\rm \rom{1}} + {\rm \rom{2}}+ {\rm \rom{3}} \right).  
\end{align*} 
To estimate ${\rm \rom{1}}$, we use $\vert P_s^D f(x) \vert \le \Vert f \Vert_{L^\infty(D)}$  so that
\begin{equation} \label{e:R1}
\begin{split}    
{\rm \rom{1}} &=\int_0^{V(|h|)V(r)} \sup_{x\in B/2} \frac{|P_s^D f(x+h)-P_s^D f(x)|}{V(|h|)} \, ds \\
&\le \int_0^{V(|h|)V(r)} \frac{2\Ve f \Ve_{L^\infty(D)}}{V(|h|)} \, ds \le c_1V(r) \Ve f \Ve_{L^\infty(D)}.
\end{split}    
\end{equation}
To estimate ${\rm \rom{2}}$, we will use Theorem \ref{t:est-p}(a). Since $s\le V(r)^2$ and $x \in B/2$, we obtain
$$  \frac{1}{d_D(x) \land 1} \lor \frac{1}{V^{-1}(\sqrt{s})} \le \frac{c_2}{V^{-1}(\sqrt{s})}. $$
Therefore, for $s \le V(r)^2$ we have 
$$|\nabla_x P_s^D f(x)| \le c_3 \left(  \frac{1}{d_D(x) \land 1} \lor \frac{1}{V^{-1}(\sqrt{s})}\right) \Ve P_s^D f\Ve_{L^\infty(D)} \le \frac{c_2c_3 }{V^{-1}(\sqrt s)} \Ve f \Ve_{L^\infty(D)} $$
for every $x \in D$. Here we used Theorem \ref{t:est-p}(a) for the first inequality. Using above inequality we conclude
\begin{equation} \label{e:R2-1}
\begin{split}
{\rm \rom{2}} &= \int_{V(|h|)V(r)}^{V(r)^2} \sup_{x\in B/2} \frac{|P_s^D f(x+h)-P_s^D f(x)|}{V(|h|)} \, ds \\
&\le \frac{|h|}{V(|h|)} \int_{V(|h|)V(r)}^{V(r)^2} \sup_{x\in B/2} |\nabla_x{P_s^D f(x^*)}| \, ds \\
&\le c_2c_3 \Ve f \Ve_{L^\infty(D)} \frac{|h|}{V(|h|)} \int_{V(|h|)V(r)}^{V(r)^2}\frac{1}{V^{-1}(\sqrt s)} \, ds,
\end{split}
\end{equation}
where $x^*$ is a point on the segment between $x$ and $x+h$. Using change of variables with $s=V^2(t)$ in the first equality and Lemma \ref{l:V} for the second inequality, we get 
\begin{align}\label{e:R2-2}
\int_{V(r)V(|h|)}^{V(r)^2} \frac{1}{V^{-1}(\sqrt s)} ds 
= 2\int_{V^{-1}(V(r)^{1/2}V(|h|)^{1/2})}^r \frac{V(t)V'(t)}{t} dt 
\le c_4 \int_\eps^r \frac{V(t)}{t} \frac{V(t)}{t} dt, 
\end{align} 
where $\eps := V^{-1}(V(|h|)^{1/2}V(r)^{1/2})$. Also, by \eqref{e:V-wsc} we have
\begin{align*}
\frac{V(t)}{V(\eps)} \le c_5\left(\frac{t}{\eps} \right)^{\alpha_2} \le c_5 \frac{t}{\eps}, \quad t \ge \eps
\end{align*}
and
$$\int_0^r \frac{V(t)}{t} dt = \int_0^r \frac{V(r)}{t} \frac{V(t)}{V(r)} dt \le c_6V(r) \int_0^r \frac{1}{t} \left(\frac{t}{r} \right)^{\alpha_1} dt \le c_7V(r). $$
Using above two inequalities, we deduce from \eqref{e:R2-2} that
\begin{align}
\begin{split} \label{e:R2-3}
\int_{V(r)V(|h|)}^{V(r)^2} \frac{1}{V^{-1}(\sqrt s)} ds  &\le c_4\int_\eps^r \frac{V(t)}{t} \frac{V(t)}{t} dt
\le c_8 \frac{V(\eps)}{\eps} \int_0^r \frac{V(t)}{t}dt \\ &\le c_9 V(r) \frac{V(\eps)}{\eps} =  c_9 V(r) \frac{V(|h|)^{1/2}V(r)^{1/2}}{V^{-1}(V(|h|)^{1/2}V(r)^{1/2})}. 
\end{split}
\end{align}
Combining \eqref{e:R2-1} and \eqref{e:R2-3}, we conclude that
\begin{align*}
&{\rm \rom{2}} \le c_{10} \Ve f \Ve_{L^\infty(D)} \frac{|h|}{V(|h|)} \cdot V(r) \frac{V(|h|)^{1/2}V(r)^{1/2}}{V^{-1}(V(|h|)^{1/2}V(r)^{1/2})} \\
&= c_{10} \Ve f \Ve_{L^\infty(D)} V(r) \frac{V(r)}{u}\frac{V^{-1}(u^2/V(r))}{ V^{-1}(u)} \le c_{11} \Ve f \Ve_{L^\infty(D)}V(r) \left( \frac{u}{V(r)} \right)^{\frac{1}{\alpha_2} -1} \le c_{11} V(r) \Ve f \Ve_{L^\infty(D)},
\end{align*} 
where $u := V(h)^{1/2}V(r)^{1/2} \le V(r)$. Here we used \eqref{e:V-inv-wsc} and $\alpha_2<1$ for the second line.

For \rom{3}, first note that for any $V(r)^2 \le s \le 1$,
$$\frac{1}{d_D(x) \land 1} \lor \frac{1}{V^{-1}(\sqrt s)} \lor \frac{1}{V^{-1}(1)} \le \frac{1}{r} \lor \frac{1}{V^{-1}(\sqrt s)} \lor \frac{1}{V^{-1}(1)} \le \frac{1}{r}.$$
So, by Theorem \ref{t:est-p}(a) we have for $V(r)^2 \leq s \leq 1$, 
\begin{equation} \label{e:R3-1}
\begin{split}
|\nabla_x p_D(s,x,y)| &\le \frac{c_{12}}{r} p_D(s,x,y) 
\le \frac{c_{13}}{r}\left(1 \land {\frac{V(d_D(x))}{s^{1/2}}}\right) \left(1 \land {\frac{V(d_D(y))}{s^{1/2}}}\right) p(s, |x-y|/4) \\
&\le \frac{c_{14}}{r} \frac{V(r)}{\sqrt s} p(s,|x-y|/4).
\end{split}
\end{equation}
Here in the second line we used $V(d_D(x)) \le c_{15}V(r)$, which follows from \eqref{e:V-wsc} and $d_D(x) \le 2r$. Thus, we obtain
\begin{equation} \label{e:1}
\begin{split}
|P^D_sf(x+h)&-P^D_sf(x)| = |h| |\nabla_x P^D_s f(x_*)|
\le |h| \Ve f \Ve_{L^\infty(D)} \int_{D} |\nabla_x p_D(s,x^*,y) | \, dy \\
&\le c_{16} \ve h \ve \Ve f \Ve_{L^\infty(D)} \frac{V(r)}{r \sqrt s} \int_{D} p\left(s,\frac{|x^*-y|}{4} \right) \, dy \le c_{17} \ve h \ve \Ve f \Ve_{L^\infty(D)} \frac{V(r)}{r \sqrt s},
\end{split}
\end{equation}
where $x_*$ is a point on the line segment between $x$ and $x+h$. Here we used $\int_{\Rn}p(s,y/4)dy =4^n$ for the last inequality. 

For $s \ge 1$, using Theorem \ref{t:est-p}(b) we have
\begin{equation}\label{e:R3-2}
|\nabla_x p_D(s,x,y)| \le \frac{c_{18}}{r} p_D (s,x,y) \le \frac{c_{19}}{r} e^{-\lambda_1 s} V(d_D(x))V(d_D(y)) \le \frac{c_{20}V(r)}{r} e^{-\lambda_1 s},
\end{equation}
Here we used $d_D(x) \le 2r$, $d_D(y) \le 1$ and \eqref{e:V-wsc} in the last inequality. Thus we arrive
\begin{equation} \label{e:2}
\begin{split}
|P^D_sf(x+h)-P^D_sf(x)| &= |h| |\nabla_x P^D_s f(x_*)|
\le |h| \Ve f \Ve_{L^\infty(D)} \int_{D} |\nabla_x p_D(s,x^*,y) | \, dy \\
&\le c_{21} \ve h \ve \Ve f \Ve_{L^\infty(D)} \frac{V(r)}{r} \int_D e^{-\lambda_1 s} \, dy \le c_{22} \ve h \ve \Ve f \Ve_{L^\infty(D)} \frac{V(r)}{r} e^{-\lambda_1 s},
\end{split}
\end{equation}
where $x_*$ is a point on the line segment between $x$ and $x+h$.

Now combining \eqref{e:1} and \eqref{e:2}, we obtain
\begin{equation} \label{e:R3-3}
\begin{split}
{\rm \rom{3}}&=\int_{V(r)^2}^\infty \frac{|P_s^D f(x+h)-P_s^D f(x)|}{V(|h|)} \,ds= \left(\int_{V(r)^2}^1 + \int_1^\infty \right) \frac{|P_s^D f(x+h)-P_s^D f(x)|}{V(|h|)}ds \\ 
&\le c_{23} \frac{V(r)}{r}\frac{|h|}{V(|h|)} \Ve f \Ve_{L^\infty(D)} \left( \int_{V(r)^2}^1 \frac{1}{ \sqrt s} \, ds + \int_1^\infty e^{-\lambda_1 s} \, ds \right) \\
&\le c_{24} \Ve f \Ve_{L^\infty(D)}(2-2V(r)+\lambda_1^{-1}).
\end{split}    
\end{equation}
The last inequality follows from $\frac{V(r)}{V(|h|)} \le c_{25} \big( \frac{r}{|h|} \big)^{\alpha_2} \le c_{25} \big( \frac{r}{|h|} \big)$ since $|h| \le r$.

Combining \eqref{e:R1}, \eqref{e:R2-3} and \eqref{e:R3-3}, we conclude
$$ [Rf]_{C^V(B/2)} \le c_{26}(1+V(r))\Ve f \Ve_{L^\infty(D)} \le c_{26}(1+V(1)) \Ve f \Ve_{L^\infty(D)}. $$
Above inequality and that $\Ve Rf \Ve_{C^V(B/2)}= [Rf]_{C^V(B/2)} + \Ve Rf \Ve_{C(B/2)}$ finish the proof.\qed

      We next provide an upper bound of $R^D f$ near the boundary. In the proof we apply the estimates on the Green function in \cite[Theorem 1.6]{GKK}.
      
      \begin{lemma} \label{l:u}
      	There exists a constant $C=C(n,a_1,a_2,\alpha_1,\alpha_2,D,\Phi(1)) > 0$ such that
      	\begin{align*}
      	\ve R^D f(x) \ve \leq C \Ve f \Ve_{L^\infty(D)} V(\diam(D)) V(d_D(x))
      	\end{align*}
      	for any $f \in L^\infty(D)$ and $x \in D$.
      \end{lemma}
      
      \pf
      The estimate on the Green function in \cite[Theorem 1.6]{GKK} and \eqref{e:V-asymp} give that for any $x,y \in D$,
      \begin{equation} \label{e:Green}
      \begin{split}
      G^D (x, y) &\leq c_1 \frac{\varphi(\ve x-y \ve)}{\ve x-y \ve^n} \left( 1 \wedge \frac{\varphi(d_D(x))}{\varphi(\ve x-y \ve)} \right)^{1/2} \left( 1 \wedge \frac{\varphi(d_D(y))}{\varphi(\ve x-y \ve)} \right)^{1/2} \\
      &\leq c_1 \frac{\varphi(\ve x-y \ve)^{1/2}}{\ve x-y \ve^n} \varphi(d_D(x))^{1/2} \le c_2 \frac{V(\ve x-y \ve)}{\ve x-y \ve^n} V(d_D(x)).
      \end{split}
      \end{equation}
  Substituting \eqref{e:Green} to \eqref{e:R^D} we obtain
      \begin{align} 
      \ve R^D f(x) \ve &\leq c_3 \Ve f \Ve_{L^\infty(D)} V(d_D(x)) \int_{D} \frac{V(\ve x-y \ve)}{\ve x-y \ve^n} dy.
      \end{align}
     Also, using \eqref{e:V-wsc} we have
      \begin{align}\begin{split}
      \int_{D} \frac{V (\ve x-y \ve)}{\ve x-y \ve^n} \, dy &\leq \int_{B(x,\diam(D))}\frac{V (\ve x-y \ve)}{\ve x-y \ve^n} \, dy \le c_4 \int_0^{\diam(D)} \frac{V(r)}{r} \, dr \\ &\leq c_5 \frac{V(\diam(D))}{\diam(D)^{\alpha_1}} \int_0^{\diam(D)} r^{\alpha_1 - 1} \, dr \leq c_6 V(\diam(D)).
      \end{split}\end{align}
      Combining above two inequalities we have proved the lemma.
      \qed
      \begin{remark} \label{r:3.7}
      	As a corollary of Lemma \ref{l:u}, we have
      	\begin{align*}
      	\Ve R^D f \Ve_{L^\infty(D)} \leq C\Ve f \Ve_{L^\infty(D)}.
      	\end{align*}
      	Hence we can simplify \eqref{e:R} to
      	\begin{align} \label{e:CV}
      	\Ve R^D f \Ve_{C^{V}(B/2)} \leq \tilde{C} \Ve f \Ve_{L^\infty(D)}
      	\end{align}
      	for some constant $\tilde{C} = \tilde{C}(n, a_1,a_2, \alpha_1, \alpha_2, D,\Phi(1)) > 0$. 
      \end{remark}

      Now we are ready to prove Theorem \ref{t:1} for the function $R^D f$. 
      \begin{proposition}
      	\label{p:R^D}
      	Assume $f \in L^\infty(D)$. Then, $R^D f \in C^{V}(D)$ and there exists a constant $C>0$ such that
      	\begin{equation}
      	\label{e:RD}
      \Ve R^D f \Ve_{C^{V}(D)} \le C \Ve f \Ve_{L^\infty(D)}.
      	\end{equation}
      	The constant $C>0$ depends only on $n,a_1,a_2,\alpha_1,\alpha_2,D$ and $\Phi(1)$.
      \end{proposition}
      
 \pf
      By \eqref{e:CV} we have
      \begin{align} \label{e:Holder ineq}
      \ve R^Df(x) - R^Df(y) \ve &\leq c_1 \Ve f \Ve_{L^\infty(D)} V(\ve x-y \ve)
      \end{align}
      for all $x, y$ satisfying $\ve x-y \ve < d_D(x)/2$. We want to show that (\ref{e:Holder ineq}) holds, perhaps with a bigger constant, for all $x,y \in D$.
      
      Let $(R_0, \Lambda)$ be the $C^{1,1}$ characteristics of $D$. Then $D$ can be covered by finitely many balls of the form $B(z_i, d_D(z_i)/2)$ with $z_i \in D$ and finitely many sets of the form $B(z^*_j, R_0) \cap D$ with $z^*_j \in \partial D$. Thus, it is enough to show that \eqref{e:Holder ineq} holds for all $x, y \in B(z_j^*, R_0) \cap D$ possibly with a larger constant.
      
      Fix $B(z_0^*, R_0) \cap D$ and assume that the outward normal vector at $z_0$ is $(0, \cdots, 0, -1)$. This is possible because the operator is invariant under the rotation. Now let $x=(x', x_n)$ and $y = (y', y_n)$ be two points in $B(z_0^*, R_0) \cap D$, and let $r = \ve x - y \ve$. Let us define for $k \geq 0$
      \begin{align*}
      x^k = (x', x_n + \lambda^k r) \quad \text{and} \quad y^k = (y', y_n + \lambda^k r),
      \end{align*}
      for some $1- 2^{-1} (1+\Lambda^2)^{-1/2} \leq \lambda <1$. Since $(1+\Lambda^2)^{-1/2} (x^k)_n \leq d_D(x^k)$, we have
      \begin{align*}
      \ve x^k - x^{k+1} \ve = \lambda^k (1-\lambda) r \leq \frac{1}{2\sqrt{1+\Lambda^2}} (x^k)_n \leq \frac{1}{2} d_D(x^k).
      \end{align*}
      Thus, we have from \eqref{e:Holder ineq} that
      \begin{align*}
      \ve R^Df(x^k) - R^Df(x^{k+1}) \ve \leq c_1 \Ve f \Ve_{L^\infty(D)} V ( \ve x^k - x^{k+1} \ve ) = c_1 \Ve f \Ve_{L^\infty(D)} V ( \lambda^k (1- \lambda) r )
      \end{align*}
      and similarly that $\ve R^Df(y^k) - R^Df(y^{k+1}) \ve \leq c_1 \Ve f \Ve_{L^\infty(D)} V ( \lambda^k (1- \lambda) r )$. Moreover, note that the distance from the line segment joining $x^0$ and $y^0$ to the boundary $\partial D$ is more than $r(1-\Lambda/2)$. Thus, this line can be split into finitely many line segments of length less than $r(1-\Lambda/2)/2$. The number of small line segments depends only on $\Lambda$. Therefore, we have $\ve R^Df(x^0) - R^Df(y^0) \ve \leq c_2 \Ve f \Ve_{L^\infty(D)} V(r)$ and hence
      \begin{align*}
      &\ve R^Df(x) - R^Df(y) \ve \\ &\leq \ve R^Df(x^0) - R^Df(y^0) \ve + \sum_{k \geq 0} \big( \ve R^Df(x^k) - R^Df(x^{k+1}) \ve + \ve R^Df(y^k) - R^Df(y^{k+1}) \ve \big) \\
      &\leq c_3 \Ve f \Ve_{L^\infty(D)} \big( V(r) + \sum_{k \geq 0} V( \lambda^k(1- \lambda) r) \big) \\
      &\leq c_4 \Ve f \Ve_{L^\infty(D)} V(r) \bigg( 1 + c_5 \sum_{k \geq 0} \big( \lambda^k (1-\lambda) \big)^{\alpha_1} \bigg) \\
      &\leq c_6 \Ve f \Ve_{L^\infty(D)} V(r).
      \end{align*}
      Recall that $r=|x-y|$. This finishes the proof. \qed
      
      In the next subsection, we will prove that the function $u=-R^D f$ is the unique viscosity solution for \eqref{e:pde1} when $f \in C(D)$.

     \subsection{Nonlocal operator and infinitesimal generator} \label{s:MH}
     
     In this section we establish the relation between viscosity solutions of \eqref{e:pde1} and solutions of the following:
     \begin{equation}
     \label{e:pde3}
     \begin{cases}
     Au=f &\text{in} ~ D , \\
     u=0 &\text{in} ~ \Rn \backslash D.
     \end{cases}
     \end{equation}
     In \cite{BLM}, the authors discussed the relation between operators $A$ and $L$, for instance, domain or values of the operators; see \cite{BLM} for the application to heat equations.
     
     At the beginning of this section we apply the strategies in \cite{BLM} to our settings and obtain some related properties. After then, we obtain comparison principle for the viscosity solution. Combining these results, we finally obtain the existence and uniqueness for Dirichlet problems \eqref{e:pde1} and \eqref{e:pde3}. Moreover, these two solutions coincide under some conditions. Also, in Section \ref{s:Harnack} we obtain Harnack inequality, which is one of the key ingredients for the standard argument of Krylov in \cite{Kry}. In Section \ref{s:pf of Thm 1.2} we will make use of Harnack inequality and the comparison principle to prove Theorem \ref{t:2}.  \\
     
     Let $D \subset \R^n$ be a bounded $C^{1,1}$ open set and let
     $$\DD=\DD(D):= \{ u \in C_0(D) : Au \in C(D)   \} $$
   be the domain of operator $A$. Recall that by \cite[Lemma 2.6]{BLM} we have
     \begin{equation}
     \label{e:gen2} Au(x) = Lu(x)
     \end{equation}
     for any $u \in C^2(x) \cap C_0(\Rn)$, $x \in D$. We first show that $u=-R^D f$ satisfies \eqref{e:pde3} when $f$ is continuous.
    
     \begin{lemma}\label{t:A} Let $f  \in C(D)$ and define $u = - R^D f$. Then, $u$ is a solution for \eqref{e:pde3}.
     \end{lemma}
     
     \pf First we claim that for any $u \in C_0(D)$ and $x \in D$,
     \begin{equation}
     \label{e:A}
     Au(x) = \lim_{t \downarrow 0} \frac{P_t^D u(x) - u(x)}{t}.
     \end{equation}
     To show \eqref{e:A}, we follow the proof in \cite[Theorem 2.3]{BLM}. Note that our domain of operator is slightly different from it in \cite[(2.8)]{BLM}. 
     
    We first observe that for any $u \in \DD$ and $x \in D$,
    \begin{align*}
    P_t^D u(x) - P_t u(x) &= \E^x u(X_t^D) - \E^x u(X_t) \\
    &= \E^x [u(X_t^D) \1_{\{\tau_D \ge t \}}] -\E^x [u(X_t) \1_{\{\tau_D \ge t \}}] - \E^x [u(X_t) \1_{\{\tau_D < t \}}] \\ &= - \E^x [u(X_t)\1_{\{\tau_D < t\}}].
    \end{align*}
    Indeed, the first and the third term in the second line cancel. Hence
    \begin{equation}
    \label{e:P}
  \frac{P_t^D u(x)-u(x)}{t} - \frac{P_t u(x)- u(x)}{t} = -\frac{\E^x [ u(X_t) \1_{\{\tau_D <t \}}]}{t}= \frac{\E^x[\big(u(X_{\tau_D}) - u(X_t)\big) \1_{\{\tau_D <t\}}]}{t}. 
  \end{equation}
  Meanwhile, by the strong Markov property we obtain
    \begin{align*}
   \left| \E^x\left[ \big(u(X_{\tau_D}) - u(X_t)\big) \1_{\{\tau_D <t\}}\right] \right| \le  \E^x\left[ \left|\E^{X_{\tau_D}}[u(X_0)-u(X_{t-\tau_D})] \right| \1_{\{\tau_D <t \}}  \right].
    \end{align*}
    Since $u \in C_0(D)$ is uniformly continuous, with stochastic continuity of L\'evy process we have that for any $\eps>0$ there is $\delta=\delta(\eps)>0$ such that
    $$ | \E^z[u(X_s)]-u(z)| < \eps  $$
    for any $z \in D$ and $0<s \le \delta$. Combining above two equations we conclude
    \begin{align*}
     \big| \E^x[\big(u(X_{\tau_D}) - u(X_t)\big) \1_{\{\tau_D <t\}}] \big| \le \eps \P^x (\tau_D < t) 
    \end{align*}
    for $0<t \le \delta$. Since $D$ is open, for any $x \in D$ we have a constant $r_x>0$ such that $B(x,r_x) \subset D$. Using \cite[Theroem 5.1 and Proposition 2.27(d)]{BSW} there exists some $M>0$ such that
    $$\frac{\P^x(\tau_D < t)}{t} \le \frac{\P^x (\tau_{B(x,r_x)} < t)}{t} \le M \quad \mbox{for all} \quad t>0.$$
Combining above inequalities we obtain that
    \begin{align*}
     \lim_{t \downarrow 0} \left| \frac{P_t^D u(x) - u(x)}{t} - Au(x) \right| &= \lim_{t \downarrow 0} \left| \frac{P_t^D u(x)-u(x)}{t} - \frac{P_t u(x)- u(x)}{t} \right| \\
     &\le  \eps \lim_{t \downarrow 0}\frac{\P^x [\tau_D < t]}{t} \le \eps M.
     \end{align*}
     Since $\eps>0$ is arbitrarily, this concludes the claim.  
     
     Now we prove the lemma. Note that $u=0$ in $D^c$ immediately follows from the definition of $R^D$. Then, by \eqref{e:A} and \eqref{e:R^D} we have that for $x \in D$,
     \begin{align}
     \begin{split}\label{e:Au}
     Au(x) &= A(-R^Df)(x) = -\lim_{t \downarrow 0}\frac{ P^D_t (R^D f) (x) - R^D f(x)}{s} \\
     &=-\lim_{t \downarrow 0} \frac{1}{t} \left[ P^D_t \Big( \int_0^\infty P_s^D f(\cdot)ds \Big)(x) -   \int_0^\infty P_s^D f(x)ds \right] \\
     &=\lim_{t \downarrow 0} \frac{1}{t} \left(   -\int_0^\infty P_{t+s}^D f(x)ds  +   \int_0^\infty P_s^D f(x )ds  \right) \\
     &= \lim_{t \downarrow 0} \frac{1}{t} \left(  - \int_t^\infty P_{s}^D f(x)ds  +  \int_0^\infty P_s^D f(x )ds  \right) \\
     &= \lim_{t \downarrow 0} \frac{\int_0^t P_s^D f(x)ds}{t} = f(x).
     \end{split}
     \end{align}
     Indeed, the third line follows from the semigroup property $P_s^D P_t^D= P_{s+t}^D$ and that $R^D f \in C_0(D)$ which follows from Proposition \ref{p:R}. This finishes the proof. \qed
     
     The next lemma shows that every solution of \eqref{e:pde3} is a viscosity solution of \eqref{e:pde1}.
     \begin{lemma}
     	\label{l:rel}
     	Assume that $f \in C(D)$ and $u \in \DD$ satisfies $Au = f$ in $D$. Then, $u$ is a viscosity solution of $Lu = f$.
     \end{lemma}
     \pf For any $x_0 \in D$ and test function $v \in C^2(\R^n)$ with $v(x_0) = u(x_0)$ and $v(y) > u(y)$ for $y \in \R^n \setminus \{x_0\}$, we have
     $$ Av(x_0) = Lv(x_0). $$
Since $v(x_0) = u(x_0)$ and $P_t^D v(x_0) \ge P_t^D u(x_0)$ for every $t>0$, we have
     $$ Av(x_0) = \lim_{t \downarrow 0} \frac{P_t^D v(x_0) - v(x_0)}{t} \ge \lim_{t \downarrow 0} \frac{P_t^D u(x_0) - u(x_0)}{t} = Au(x_0). $$
     Thus, we arrive
     $$ Lv(x_0) \ge Au(x_0), $$
     which concludes that $u$ is a viscosity solution of \eqref{e:pde1}. \qed
    
     Now we see comparison principle in \cite{CS2}. This implies the uniqueness of viscosity solution for \eqref{e:pde1}.

     \begin{theorem} [Comparison principle] \label{t:cp}
     	Let $D$ be a bounded open set in $\Rn$. Let $u$ and $v$ be bounded functions satisfying $Lu \geq f$ and $Lv \leq f$ in $D$ in viscosity sense for some continuous function $f$, and let $u \leq v$ in $\Rn \setminus D$. Then $u \leq v$ in $D$.
     \end{theorem}
     
     \pf We first claim that $L$ satisfies \cite[Assumption 5.1]{CS2}. More precisely, there exists constant $r_0 \geq 1$ such that for every $r \geq r_0$, there exists a constant $\delta= \delta(r) > 0$ satisfying $Lw > \delta$ in $B_r$, where $w(x) = 1 \land \frac{ \ve x \ve^2}{r^3}$.

     Let $r_0=4$, $r \geq 4$ and $x \in B_r$. Note that by $r \ge 4$ we have
     $$\frac{|y|^2}{r^3} \le \frac{4r^2}{r^3} \le 1, \qquad  y \in B_{2r}. $$
     Thus, for $y \in B_r$ we obtain
     \begin{align*}
     w(x+y) + w(x-y) - 2w(x) = \frac{\ve x +y \ve^2 + \ve x - y \ve^2 - 2\ve x \ve^2}{r^3} = \frac{2\ve y \ve^2}{r^3}.
     \end{align*}
       On the other hand, for $y \in B_r^c$ we have $$w(x+y) + w(x-y) - 2w(x) \geq \frac{2|y|^2}{r^3} \land (1 - 2w(x)) > 0.$$
     Therefore, since $w \in C^2(\R^d)$ we have
     \begin{align*}
    &Lw(x):= \frac{1}{2}\int_{\R^n} \left( w(x+y) + w(x-y) - 2w(x) \right) J(y) \, dy \\
    &= \frac{1}{2}\int_{B_r} \left( w(x+y) + w(x-y) - 2w(x) \right) J(y)\, dy + \frac{1}{2}\int_{B_r^c} \left( w(x+y) + w(x-y) - 2w(x) \right) J(y) \, dy
     \\ &\ge  \frac{1}{r^3} \int_{B_r} |y|^2 J(y)dy \, =: \delta(r)>0
     \end{align*}
     for every $r \geq r_0=4$ and $x \in B_r$. Since $L$ satisfies \cite[Assumption 5.1]{CS2}, we can apply Theorem 5.2 therein, which proves the theorem.
     \qed
     
           The following uniqueness of viscosity solution is immediate.
           
           \begin{corollary} \label{c:uniqueness}
           	Let $D$ be a bounded open set in $\Rn$ and let $f \in C(D)$. Then there is at most one viscosity solution of \eqref{e:pde1}.
           \end{corollary}     
                 Here is the main result in this section.
           \begin{theorem}\label{t:sol}
           	Assume that $f \in C(D)$. Then, $u=-R^D f \in \DD$ is the unique solution of \eqref{e:pde3}. Also, $u$ is the unique viscosity solution of \eqref{e:pde1}.
           	\end{theorem}
           	\pf By Lemma \ref{t:A}, we have that $u=-R^D f \in \DD$ is solution of \eqref{e:pde3}. Now, Lemma \ref{l:rel} and Corollary \ref{c:uniqueness} conclude the proof. \qed

          \noindent  \textbf{Proof of Theorem \ref{t:1}} By Theorem \ref{t:sol}, the unique viscosity solution for \eqref{e:pde1} is given by $u=-R^Df$. Therefore, Proposition \ref{p:R^D} yields the H\"older regularity of viscosity solution with respect to $C^V$-norm. By \eqref{e:V-asymp}, we have $V \asymp \overline{\phi}$ and this concludes the proof.  \qed

\section{Boundary regularity} \label{s:bdry reg}

\subsection{Barriers} \label{s:barriers}

Throughout this section, $D \subset \Rn$ is a bounded $C^{1,1}$ open set. Without loss of generality, we assume that $\diam(D) \leq 1$.
Since $d_D$ is only $C^{1,1}$ near $\partial D$, we need to consider the following ``regularized version" of $d_D$.

\begin{definition}\label{d:psi}
We call $\psi : D \rightarrow (0, \infty)$ the regularized version of $d_D$ if $\psi \in C^{1,1}(D)$ and it satisfies 
	\begin{equation}\label{e:psi1}
\tilde{C}^{-1}d_D(x) \le \psi(x) \le \tilde{C}d_D(x), \quad \Ve\nabla \psi(x)\Ve \le \tilde{C} \quad \mbox{and} \quad \Ve \nabla \psi(x) - \nabla \psi(y)\Ve \le \tilde{C}|x-y|
	\end{equation}
	for any $x,y \in D$, where the constant $\tilde{C}>0$ depends only on $D$. 
\end{definition}
 For $D=B(0,1)$, there exists a regularized version of $d_{B(0,1)}$ which is $C^2$ and isotropic. Denote this function by $\Psi$ and let $C=C(n)$ be the constant in \eqref{e:psi1} for the function $\Psi$. For any open ball $B_r:=B(x_0,r)$, we will take the regularized version of $d_{B_r}$ which is defined by
$\Psi_r(x):= \Psi(\frac{x-x_0}{r})$. Then, $\Psi_r$ satisfies 
\begin{equation}\label{e:psi2}
C^{-1}d_{B_r}(x) \le \Psi_r(x) \le Cd_{B_r}(x), \quad \Ve \nabla \Psi_r\Ve \le C \quad \mbox{and} \quad \Ve \nabla^2\Psi_r(x)\Ve \le \frac{C}{r}
\end{equation}
for any $x,y \in B(x_0,r)$. The last estimate follows from the fact that $\Psi \in C^2(B_r)$.

We first introduce the following three lemmas which will be used to construct a barrier for $L$.

\begin{lemma}\label{l:3.4}
Assume that $D$ is a bounded $C^{1,1}$ open set and let $\psi$ be a regularized version of $d_D$. Then, for every $x \in \Rn$ and $x_0 \in D$ we have
	\begin{equation}\label{e:3.4}
|\psi(x)-(\psi(x_0)+\nabla \psi(x_0) \cdot (x-x_0))_+| \le \tilde{C}|x-x_0|^2
\end{equation}
	where $\tilde{C}$ is the constant in \eqref{e:psi1}. 
	In addition, when $D=B(0,r)$ and $\psi=\Psi_r$ we have \eqref{e:3.4} with $\tilde{C}=\frac{C}{r}$ where $C$ is the constant in \eqref{e:psi2}.
\end{lemma}

\pf Let $\tilde{\psi}$ be a $C^{1,1}$ extension of $\psi |_D$ satisfying $\tilde{\psi} \le 0$ in $\R^n \backslash D.$ Then, since $\tilde{\psi} \in C^{1,1}(\R^n)$ we clearly have
\begin{equation}\label{e:4.3}
|\tilde{\psi} (x) - \psi(x_0) - \nabla \psi(x_0) \cdot (x-x_0)|=|\tilde{\psi} (x) - \tilde{\psi}(x_0) - \nabla \tilde{\psi} \cdot (x-x_0)| \le \tilde{C}|x-x_0|^{2} 
\end{equation}
in all of $x \in \R^n$. Using $|a_+-b_+| \le |a-b|$ and $ (\tilde{\psi})_+=\psi$, we have
$$
|\psi(x)-(\psi(x_0)+\nabla \psi(x_0) \cdot (x-x_0))_+| \le |\tilde{\psi} (x) - \psi(x_0) - \nabla \psi(x_0) \cdot (x-x_0)| \\ 
\le \tilde{C}|x-x_0|^2
$$
for all $x\in \R^n$. If $D=B(0,r)$ and $\psi=\Psi_r$, the constant $\tilde{C}$ in \eqref{e:4.3} become $\frac{C}{r}$. Thus, the conclusion of lemma follows. \qed

Next lemma is a collection of inequalities which will be used for this section. Note that we can easily check these inequalities when $\varphi(r)=r^{2\alpha}$ and $V(r)=r^{\alpha}$ with $0<\alpha<1$. The inequalities \eqref{e:varphi-inf} and \eqref{e:V-inf} are in \cite[Lemma 3.5]{BGR2}. We provide the proof for the completeness. 
\begin{lemma}
	\label{l:phi-inf}
	There exists a constant $C_1=C_1(n, a_1, \alpha_1,\alpha_2)>0$ such that for any $0<r \le 1$,
	\begin{equation}\label{e:varphi-0}
	\int_0^r \frac{s}{\varphi(s)}ds \le \frac{C_1 r^2}{\varphi(r)},
	\end{equation}
	\begin{equation}
	\label{e:varphi-inf}
	\int_{r}^\infty \frac{1}{s\varphi(s)}ds \le \frac{C_1}{\varphi(r)},
	\end{equation}
		\begin{equation}\label{e:V-0}
		\int_0^r \frac{1}{V(s)}ds \le \frac{C_1 r}{V(r)}, \quad \int_0^r \frac{V(s)}{s}ds \le C_1 V(r) 
		\end{equation}
		and
		\begin{equation}
		\label{e:V-inf}
		\int_{r}^\infty \frac{V(s)}{s\varphi(s)}ds \le \frac{C_1}{V(r)}.
		\end{equation}
\end{lemma}
\pf The inequalities \eqref{e:varphi-0} and \eqref{e:V-0} can be proved using weak scaling conditions \eqref{e:varphi-wsc} and \eqref{e:V-wsc}: by \eqref{e:varphi-wsc}, we have
$$\int_0^r \frac{s}{\varphi(s)}ds = \int_0^r \frac{s}{\varphi(r)} \frac{\varphi(r)}{\varphi(s)}ds \le c_1\int_0^r \frac{s}{\varphi(r)}  \Big(\frac{r}{s}\Big)^{2\alpha_2} ds = \frac{c_1}{2-2\alpha_2} \frac{r^2}{\varphi(r)},$$
and by \eqref{e:V-wsc} we have
$$\int_0^r \frac{1}{V(s)}ds = \int_0^r \frac{1}{V(r)} \frac{V(r)}{V(s)} ds \le \int_0^r c_2 \Big( \frac{r}{s} \Big)^{\alpha_2} ds = \frac{c_2}{1-\alpha_2} \frac{r}{V(r)} $$
and
$$\int_0^r \frac{V(s)}{s}ds = \int_0^r  \frac{V(r)}{s} \frac{V(s)}{V(r)} ds \le \int_0^r \frac{V(r)}{s} c_2 \Big( \frac{s}{r} \Big)^{\alpha_1} ds = \frac{c_2}{\alpha_1} V(r). $$
 Let $\mathcal{P}(r):= \int_\R \big( 1 \land \frac{|x|^2}{r^2} \big) J(x)dx$ be the Pruitt function of $X$. By \cite[(6) and Lemma 1]{BGR1} and \eqref{e:V-asymp}, we have a constant $c_3>0$ satisfying
\begin{equation}
\label{e1}\mathcal{P}(r) \le c\varphi(r)^{-1} \le c_3 V(r)^{-2}, \quad r>0.
\end{equation}
Let $\mathcal{P}_1(r):= \int_r^\infty \frac{1}{s\varphi(s)}ds$. Note that we have
\begin{equation}
\label{e2}
\mathcal{P}_1(r)=\omega_n^{-1} \int_{B(0,r)^c} \Big( 1 \land \frac{|x|^2}{r^2} \Big) J(|x|)dx \le \omega_n^{-1} \mathcal{P}(r) \le c_4 V(r)^{-2}, \quad r>0.
\end{equation}
Thus, \eqref{e2} and \eqref{e:V-asymp} imply \eqref{e:varphi-inf}.
Also, using integration by parts and \eqref{e2} we have
\begin{align*}
\int_r^{\infty} \frac{V(s)}{s\varphi(s)} ds &= \int_r^\infty V(s) d(-\mathcal{P}_1)(s) \\
&= V(r) \mathcal{P}_1(r) - \lim_{s \rightarrow \infty} V(s) \mathcal{P}_1(s) + \int_r^\infty V'(s) \mathcal{P}_1(s)ds \\
&\le c_5 \left(\frac{1}{V(r)} - \lim_{s \rightarrow \infty} \frac{1}{V(s)} + \int_r^\infty \frac{V'(s)}{V(s)^2}ds  \right) =  \frac{2c_5}{V(r)},
\end{align*}
which concludes \eqref{e:V-inf}.  \qed

\begin{lemma}\label{l:3.5} Let $U \subset \R^n$ be a $C^{1,1}$ open set, which can be unbounded. Then there exists a constant $C_2=C_2(n,U, a_1,a_2,\alpha_1, \alpha_2)>0$ such that for any $x \in U$ and $0<r \le 1$, 
	\begin{equation}\label{e:4.4.1}
	\int_{ U \cap \big(B(x,r)  \backslash B(x,d_U(x)/2) \big)} \frac {V(d_U(y))}{d_U(y)} \frac{dy}{|x-y|^{n-2} \varphi(|x-y|)} \le \frac{C_2 r}{V(r)} .
	\end{equation}
\end{lemma}

\pf Fix $x \in U$ and denote $\rho:=d_U(x)<2r$, $B_r := B(x,r)$ for $r > 0$ and $B_r = \emptyset$ for $r\le 0$. First note that there is a constant $\kappa= \kappa(U)>0$ such that the level set  \noindent $\{d_U \ge t \} = \{x \in U | d_U(x) \ge t \}$ is $C^{1,1}$ for any $t \in (0, \kappa]$ since $U$ is $C^{1,1}$. Without loss of generality we can assume $\kappa \le r$ because $\kappa$ can be arbitrarily small.

Since $B_R \cap \left\{ d_U \ge \kappa \right\} = \emptyset$ for every $R \le \kappa - \rho$, we have 
\begin{align*} 
	&\int_{(B_r \backslash B_{\rho/2}) \cap \left\{ d_U \ge \kappa \right\}}\frac{V(d_U(y))}{d_U(y)} \frac{dy}{|x-y|^{n-2} \varphi(|x-y|)} \\ &= \int_{(B_r \backslash B_{\max\{\rho/2,\kappa-\rho\}}) \cap \{ d_U \ge \kappa \}}\frac{V(d_U(y))}{d_U(y)} \frac{dy}{|x-y|^{n-2} \varphi(|x-y|)} \\ &\le
	\int_{(B_r \backslash B_{2\kappa /3}) \cap \{ d_U \ge \kappa \}}\frac{V(d_U(y))}{d_U(y)} \frac{dy}{|x-y|^{n-2} \varphi(|x-y|)},
\end{align*}
where the last line follows from $\rho/2 \lor (\kappa-\rho) \ge \frac{2\kappa}{3}.$ Using
$$\kappa \le d_U(y) \le r+\kappa \le 2r \quad \mbox{and} \quad \frac{2\kappa}{3} \le |x-y| \le r $$
for every $y \in (B_r \backslash B_{2\kappa/3}) \cap \{d_U \ge \kappa \}$, we arrive that for any $x \in U$,
\begin{align}\begin{split}
&\int_{(B_r \backslash B_{2\kappa/3}) \cap \left\{ d_U \ge \kappa \right\}}\frac{V(d_D(y))}{d_D(y)} \frac{dy}{|x-y|^{n-2} \varphi(|x-y|)} \\ &\le \int_{(B_r \backslash B_{2\kappa /3}) \cap \{ d_U \ge \kappa \}}\frac{V(2r)}{\kappa} \frac{dy}{|x-y|^{n-2} \varphi(|x-y|)} 
\\ &\le c_1\frac{V(r)}{\kappa} \int_0^r \frac{s}{\varphi(s)}ds \le c_2(\kappa) \frac{r^2}{V(r)} \le c_2(\kappa) \frac{r}{V(r)}, \label{e:3.5}
\end{split}
\end{align}
where we used \eqref{e:V-asymp} and \eqref{e:varphi-0} for the second last inequality.
Thus, it suffices to estimate the integrand \eqref{e:4.4.1} in the set $(B_r \backslash B_{\rho/2}) \cap \{ 0 < d_U < \kappa \}.$

We will utilize the following estimates on Hausdorff measure in [RV15], that is, there exists a constant $c_3(U)>0$ such that that for every $x \in U$ and $t \in (0,\kappa)\,$,
\begin{equation}\label{e:3.5-3}
\sH^{n-1}(\left\{ d_U=t \right\} \cap (B_{2^{-k+1}r} \backslash B_{2^{-k}r})) \le c_3(2^{-k}r)^{n-1}
\end{equation}
which follows from the fact that the level set $\left\{d_U=t\right\}$ is $C^{1,1}$ for $t \in (0,\kappa)$.

Let us denote $C_n := B_{r 2^{-n}}$ for $n \ge 0$ and let $M \in \N$ be the natural number satisfying $2^{-M}r \le \rho/2 \le 2^{-M+1}r.$ Using $|x-y| \ge 2^{-k}r$ for every $y \in C_{k-1} \backslash C_k$ and $\varphi$ is increasing for the third line, we have
\begin{align}
&\int_{(B_r \backslash B_{\rho/2}) \cap \{ 0<d_U < \kappa \}}\frac{V(d_U(y))}{d_U(y)} \frac{dy}{|x-y|^{n-2} \varphi(|x-y|)} \nn \\ & \le \quad \sum_{k=1}^M \int_{(C_{k-1} \backslash C_k) \cap \{ 0<d_U < \kappa \}}\frac{V(d_U(y))}{d_U(y)} \frac{dy}{|x-y|^{n-2} \varphi(|x-y|)} \nn \\
& \le \quad \sum_{k=1}^M \frac{1}{(2^{-k}r)^{n-2}\varphi(2^{-k}r)} \int_{(C_{k-1} \backslash C_k) \cap \{ 0< d_U < \kappa  \}}\frac{V(d_U(y))}{d_U(y)}dy \nn \\
& = \quad \sum_{k=1}^M \frac{1}{(2^{-k}r)^{n-2}\varphi(2^{-k}r)} \int_{(C_{k-1} \backslash C_k) \cap \{ 0< d_U < \kappa \}}\frac{V(d_U(y))}{d_U(y)}|\nabla d_U(y)|dy. \nn
\end{align}
Here we used $|\nabla d_U(y)|=1$ for $y \in \{0<d_U<\kappa\}$ for the last line. (See \cite{ROV}.) \\
For any $1 \le k \le M$ and $y \in C_{k-1}$ we have $d_U(y) \le 2^{-k+1}r + \rho \le (2^{-k+1} + 2^{-M+2})r \le 6 \cdot 2^{-k}r$, which implies $C_{k-1} \subset \{d_U < 6 \cdot 2^{-k}r \}$. Thus, combining this with above inequality we have
\begin{align}\begin{split}
&\int_{(B_r \backslash B_{\rho/2}) \cap \{ 0<d_U < \kappa \}}\frac{V(d_U(y))}{d_U(y)} \frac{dy}{|x-y|^{n-2} \varphi(|x-y|)} \\  &\le \quad \sum_{k=1}^M \frac{1}{(2^{-k}r)^{n-2}\varphi(2^{-k}r)} \int_{(C_{k-1} \backslash C_k) \cap \{ 0< d_U < 6 \cdot 2^{-k}r  \}}\frac{V(d_U(y))}{d_U(y)}|\nabla d_U(y)|dy. \label{e:3.5-2}
\end{split}
\end{align}
Plugging $u(y)=d_U(y)$ and $g(y)=\frac{V(d_U(y))}{d_U(y)}$ into the following coarea formula
$$\int_{D}g(y) | \nabla u(y)| dy = \int_{-\infty}^{\infty} \left( \int_{u^{-1}(t)} g(y)d\sH_{n-1}(y) \right) dt, $$
we obtain
\begin{align}\label{e:3.5-4}
\begin{split}
&\sum_{k=1}^M \frac{1}{(2^{-k}r)^{n-2}\varphi(2^{-k}r)} \int_{(C_{k-1} \backslash C_k) \cap \{ 0< d_U < 6 \cdot 2^{-k}r  \}}\frac{V(d_U(y))}{d_U(y)}|\nabla d_U(y)|dy  \\  &= \quad \sum_{k=1}^M \frac{1}{(2^{-k}r)^{n-2}\varphi(2^{-k}r)}  \int_0^{6 \cdot 2^{-k} r} \int_{(C_{k-1} \backslash C_k) \cap \left\{ d = t \right\}} \frac{V(t)}{t} d\mathcal{H}^{n-1}(y)dt  \\ & 
\le \quad  \sum_{k=1}^M \frac{1}{(2^{-k}r)^{n-2}\varphi(2^{-k}r)} \int_0^{6\cdot 2^{-k}r}    c_3 (2^{-k}r)^{n-1}     \frac{V(t)}{t}dt  \\& = \quad  c_3 \sum_{k=1}^M \frac{ 2^{-k}r}{\varphi(2^{-k}r)}   \int_0^{6\cdot 2^{-k}r}\frac{V(t)}{t}dt \le c_4 \sum_{k=1}^M \frac{ 2^{-k}r}{\varphi(2^{-k}r)} V(6 \cdot 2^{-k}r),
\end{split}
\end{align}
where we used \eqref{e:3.5-3} for the third line and \eqref{e:V-0} for the last line. Also, by \eqref{e:V-wsc} and \eqref{e:V-asymp},
\begin{equation}\label{e:3.5-5}
\begin{split}
\sum_{k=1}^M \frac{ 2^{-k}r}{\varphi(2^{-k}r)} V(6 \cdot 2^{-k}r) &\le \sum_{k=1}^M \frac{ 2^{-k}r}{V(2^{-k}r)} = \sum_{k=1}^M  \int_{2^{-k}r}^{2^{-k+1}r} \frac{1}{V(2^{-k}r)}ds \\
&\le \int_0^r \frac{1}{V(s)}ds \le c_5 \frac{r}{V(r)},
\end{split}
\end{equation}
where in the last two inequalities we have used that $V$ is increasing and \eqref{e:V-0}. \\

\noindent Using \eqref{e:3.5-2}, \eqref{e:3.5-4}, and \eqref{e:3.5-5}, we conclude
$$ \int_{(B_r \backslash B_{\rho/2}) \cap \{ d < \kappa \}}\frac{V(d_U(y))}{d_U(y)} \frac{dy}{|x-y|^{n-2} \varphi(|x-y|)} \le \frac{c_4c_5 r}{V(r)}.  $$
This and \eqref{e:3.5} finish the proof.   \qed

\noindent Now we are ready to show that $V(\psi)$ acts as a barrier of $L$ on $D$.
\begin{proposition}\label{p:barrier1}
	Let $L$ be given by \eqref{d:L} and $\psi$ be a regularlized version of $d_D$. Then there exists a constant $\tilde{C}_3=\tilde{C}_3(n,a_1,a_2,\alpha_1,\alpha_2,D)>0$ such that
	\begin{equation}\label{e:5.1}
	|L(V(\psi))| \le \tilde{C}_3 \quad in \,\, D.
	\end{equation}
	where $V$ is the renewal function with respect to $\Phi$. In addition, if $D=B(0,r)$ is a ball with radius r, there exists a constant $C_3=C_3(n,a_1,a_2,\alpha_1,\alpha_2)>0$ such that 
	\begin{equation}\label{e:5.2}
		|L(V(\psi))| \le \frac{C_3}{V(r)} \quad in \,\, B(0,r),
	\end{equation}
	where $\psi=\Psi_r$ is a regularized version of $d_{B(0,r)}$ defined in \eqref{e:psi2}. Note that $C_3$ is independent of r.
\end{proposition}

\pf We prove \eqref{e:5.2} only. The proof of \eqref{e:5.1} is similar. 

Let $x_0 \in B_r:=B(0,r)$ and $\rho := d_{B_r}(x_0)$. 
First we prove \eqref{e:5.2} for the case $\rho \ge \kappa r > 0$ with $\kappa=1/(8C^2)$. In this case, we have
\begin{equation} \label{e:5.3}
\begin{split}
|L(V(\psi))(x_0)|
=& \lv\int_{\Rn}\left(\frac{V(\psi(x_0+y))+V(\psi(x_0-y))}{2}-V(\psi(x_0)) \right)\frac{J(1)}{|y|^n \varphi(|y|)}dy \rv \\
\le& \int_{B_{\kappa r/2}}\lV \nabla^2 [V(\psi(x_*))]\rV \frac{J(1)}{|y|^{n-2} \varphi(|y|)}dy \\
&+\int_{B^c_{\kappa r/2}} \lv \frac{V(\psi(x_0+y))+V(\psi(x_0-y))}{2}-V(\psi(x_0)) \rv\frac{J(1)}{|y|^n \varphi(|y|)}dy, 
\end{split}
\end{equation}
where $x^*$ is a point on the segment between $x_0-y$ and $x_0+y$, so that $d_{B_r}(x_*) \ge \kappa r/2$ when $y \in B_{\kappa r /2 }$. Using \eqref{e:V-wsc}, \eqref{e:psi2}, and Lemma \ref{l:V}, we have 
$$\Ve \nabla^2 [V(\psi(x_*))] \Ve  \le |V''(\psi(x))| \Ve \nabla \psi(x) \Ve ^2+|V'(\psi(x))| \Ve \nabla^2 \psi(x) \Ve  \le \frac{c_1(\kappa) V(r)}{r^2},$$ 
which yields to estimate the first term of \eqref{e:5.3} by
\begin{align*}
 \int_{B_{\kappa r/2}}\Ve \nabla^2 [V(\psi(x_*))]\Ve \frac{J(1)}{|y|^{n-2} \varphi(|y|)}dy 
 &\le c_1 \frac{V(r)}{r^2} \int_{B_{\kappa r/2}}\frac{1}{|y|^{n-2} \varphi(|y|)}dy \\
 &= c_2 \frac{V(r)}{r^2} \int_0^{\kappa r /2 } \frac{s}{\varphi(s)} ds \le \frac{c_3}{V(r)}.
 \end{align*}
In the last inequality above, we have used \eqref{e:varphi-0}, \eqref{e:varphi-wsc}, and \eqref{e:V-asymp}. For the second term, using $\psi(x) \le C d_{B_r}(x) \le Cr$ for any $x \in B_r$, we have $$\lv \frac{V(\psi(x_0+y))+V(\psi(x_0-y))}{2}-V(\psi(x_0)) \rv \le 2 V(C r) \le c_4 V(r).$$ 
Therefore,
$$ \int_{B^c_{\kappa r/2}} \lv \frac{V(\psi(x_0+y))+V(\psi(x_0-y))}{2}-V(\psi(x_0)) \rv\frac{J(1)}{|y|^n \varphi(|y|)}dy \le c_5V(r) \int_{\kappa r /2}^{\infty} \frac{1}{ s\varphi(s)} ds \le \frac{c_6(\kappa)}{V(r)}.$$ 
In the last inequality we have used \eqref{e:varphi-inf}, \eqref{e:varphi-wsc}, and \eqref{e:V-asymp}. 
Therefore, \eqref{e:5.2} for the case $\rho \ge \kappa r$ holds with $C_3 = c_3 + c_6$. 

Now it suffices to consider the case $\rho < \kappa r$. Denote
$$l(x):=(\psi(x_0)+\nabla \psi(x_0) \cdot (x-x_0))_+,$$
which satisfies
$$ L(V(l))=0 \quad \mbox{on} \quad \{l>0\} $$
by \eqref{e:pde-V}. Note that
$ \psi(x_0)=l(x_0)$ and $\nabla \psi(x_0)=\nabla l(x_0).$
Moreover, by \eqref{e:3.4} we have
\begin{equation}\label{e:3.6.1}
|\psi(x)-l(x)| \le \frac{C}{r}|x-x_0|^{2}.
\end{equation}
For any $0<a\le b \le C$, there exists $a_* \in [a,b]$ satisfying $|V(a)-V(b)|=|a-b|V'(a_*)$. Using Lemma \ref{l:V} in the first inequality we have
$$
|V(a)-V(b)| = |a-b|V'(a_*) \le c_7|a-b| \frac{V(a_*)}{a_*} 
\le c_8 |a-b| \frac{V(a)}{a} .
$$
Here we used \eqref{e:V-wsc} with $c=C$ for the second inequality. Therefore, for any $a,b \in (0,C]$ we have
$$|V(a)-V(b)| \le c_8|a-b| \left(\frac{V(a)}{a}+\frac{V(b)}{b}\right). $$
Also, one can easily see the following inequality
\begin{equation}\label{e:3.6.2}
|V(a)-V(b)| \le c_8|a-b| \left(\frac{V(a)}{a} {\bf 1}_{\{a>0\}}+\frac{V(b)}{b}\cdot {\bf 1}_{\{b>0\}}\right)
\end{equation}
for any $0 \le a,b \le C$ by using Lemma \ref{l:V}. 
 
By \eqref{e:3.6.1} and \eqref{e:3.6.2} we have that for any $x \in B_r(x_0)$,
\begin{align}\label{e:3.6.3}
|V(\psi(x))-V(\ell(x))| &\le \frac{c_8}{r}|x-x_0|^{2}\left(\frac{V(\psi(x))}{\psi(x)} \1_{\{\psi(x)>0  \}}+\frac{V(\ell(x))}{\ell(x)} \1_{\{\ell(x)>0\}}\right) \\ &\le \frac{c_9}{r}|x-x_0|^{2}\left(\frac{V(d_{B_r}(x))}{d_{B_r}(x)} \1_{\{d_{B_r}(x)>0\}}+\frac{V(\ell(x))}{\ell(x)}\1_{\{\ell(x)>0\}}\right), \nn
\end{align}
where we used $\psi(x) \le Cd_{B_r}(x) \le C$ and $\ell(x)=(\psi(x_0)+\nabla \psi(x_0) \cdot (x-x_0))_+ \le Cd_{B_r}(x_0) + Cr \le C$ for the first inequality and \eqref{e:V-wsc} for the second. \\

On the other hand, for any $x \in B_{\rho/2}(x_0)$ with $\rho \le \kappa r$ we have
$$ |\ell(x) - \psi(x)| \le \frac{C}{r}|x-x_0|^2 \le \frac{C}{r} \rho^2 \le C \kappa \rho$$
and
$$ C^{-1} \frac{\rho}{2} \le C^{-1}d_{B_r}(x) \le \psi(x). $$
Thus, using $\kappa = 1/(8C^2)$ we obtain
$$ \frac{1}{2} \psi(x) \le \ell(x) \le 2\psi(x)   \quad \mbox{for any} \quad x \in B_{\rho/2}(x_0).$$
Using $\frac{\rho}{2} \le d_{B_r}(x) \le 2\rho$, we arrive at
$$  \psi(x), \ell(x) \in [(4C)^{-1} \rho, 4C\rho]. $$ 
Therefore, there exists $y \in ((4C)^{-1} \rho, 4C \rho)$ satisfying $$\frac{V(\psi(x)) - V(\ell(x))} {\psi(x)-\ell(x)} = V'(y),$$ so using \eqref{e:3.6.1} and \eqref{e:V-diff}, we have
\begin{align}
|V(\psi(x))-V(\ell(x))|&=|\psi(x)-\ell(x)|V'(y)   \le \frac{c_{10}}{r} |x-x_0|^2  \frac{V(y)}{y} \label{e:5.4}  \\
&\le \frac{c_{11}}{r} |x-x_0|^2  \frac{V((4C)^{-1} \rho)}{(4C)^{-1}\rho} 
\le \frac{c_{12}}{r} |x-x_0|^2 \frac{V(\rho)}{\rho} \nn
\end{align}
for $x \in B_{\rho/2}(x_0)$. Here we used \eqref{e:V-diff} and \eqref{e:V-wsc} for the second line. Also, for any $x \in B^c_r(x_0)$ we have 
$$V(\ell(x)) = V(\psi(x_0)+(x-x_0)\nabla \psi(x_0)) \le V(C\rho + C|x-x_0|) \le V(2C|x-x_0|) \le c_{13} V(|x-x_0|)$$ 
and
$$V(\psi(x))\le V(Cr) \le V(C|x-x_0|) \le c_{13} V(|x-x_0|),$$ 
where we have used \eqref{e:V-wsc} and $\rho \le r \le |x-x_0|.$
 Thus we obtain
\begin{equation}\label{e:5.5}
|V(\psi)-V(\ell)|(x) \le c_{14} V(|x-x_0|)
\end{equation}
for $x \in B^c_r(x_0)$. Therefore, by taking $x=y+x_0$ for \eqref{e:3.6.3}, \eqref{e:5.4}, and \eqref{e:5.5} we have
\begin{align}
|V(\psi)-V(\ell)|(y+x_0) \le c
\begin{cases}
\frac{1}{r}\frac{V(\rho)}{\rho}|y|^2   &\mbox{for} \ y \in B_{\rho/2} \\
\frac{|y|^2}{r}\left(\frac{V(d_{B_r}(x_0+y))}{d_{B_r}(x_0+y)} {\bf 1}_{\{d_{B_r}(x_0+y)>0\}}+\frac{V(l(x_0+y))}{l(x_0+y)}{\bf 1}_{\{l(x_0+y)>0\}}\right) &\mbox{for} \ 
y \in B_r \backslash B_{\rho/2} \\
V(|y|) &\mbox{for} \ y \in B^c_r
\end{cases} \nn
\end{align}
where $c = c_9 \lor c_{12} \lor c_{14}.$
Hence, recalling that $L(V(\ell))(x_0)=0$ and $\psi(x_0)=\ell(x_0)$, we find that
\begin{align*}
|L(V(\psi))&(x_0)|= |L(V(\psi)-L(V(\ell)))(x_0)| \\
& = \int_{\Rn} |V(\psi)-V(\ell)|(x_0+y) \frac{J(1)}{|y|^n \varphi(|y|)}dy  \\
& \le \frac{c}{r} \frac{V(\rho)}{\rho} \int_{B_{\rho/2}}|y|^2  \frac{J(1)}{|y|^n \varphi(|y|)}dy + c\int_{B^c_r} V(|y|)\frac{J(1)}{|y|^n \varphi(|y|)}dy \\
&\quad + c\int_{B_r \backslash B_{\rho/2}} \frac{|y|^2}{r}\left(\frac{V(d_{B_r}(x_0+y))}{d_{B_r}(x_0+y)} {\bf 1}_{\{d_{B_r}(x_0+y)>0\}}+\frac{V(\ell(x_0+y))}{\ell(x_0+y)}{\bf 1}_{\{\ell(x_0+y)>0\}}\right)\frac{J(1)}{|y|^n \varphi(|y|)}dy \\
&=: {\rm {\rom{1}}}+{\rm {\rom{2}}}+{\rm {\rom{3}}}.
\end{align*}
For \rom{1}, using \eqref{e:varphi-0} we have
\begin{align*}
{\rm {\rom{1}}} &= \frac{c}{r} \frac{V(\rho)}{\rho} \int_{B_{\rho/2}}|y|^2  \frac{J(1)}{|y|^n \varphi(|y|)}dy = \frac{c_{15}}{r} \frac{V(\rho)}{\rho} \int_0^{\rho/2} \frac{s}{\varphi(s)}ds  \\ &\le \frac{c_{16}}{r} \frac{V(\rho)}{\rho} \frac{(\rho/2)^2}{\varphi(\rho/2)} \le \frac{c_{17}}{V(r)}  \left( \frac{\rho}{r} \frac{V(r)}{V(\rho)} \right) \le \frac{ c_{18} }{V(r)},
\end{align*}
where we used \eqref{e:V-asymp} and \eqref{e:V-wsc} for the last two inequalities. Also, using \eqref{e:V-inf} we obtain
\begin{align*}
{\rm {\rom{2}}}= c\int_{B^c_r} V(|y|)\frac{J(1)}{|y|^n \varphi(|y|)}dy = c_{19} \int_r^\infty \frac{V(s)}{s\varphi(s)} ds \le \frac{c_{20}}{V(r)}.
\end{align*}
For the estimate of \rom{3}, we first observe that for any $y \in \{ \ell >0\} := H$,
$$ \lv \frac{\ell(y)}{d_H(y)} \rv = \Ve \nabla \psi(x_0) \Ve \le C. $$
Thus, by \eqref{e:V-wsc} we have
$$ \frac{V(\ell(y))}{\ell(y)} \le c_{21} \frac{V(C d_H(y))}{d_H(y)} \le c_{22} \frac{V(d_H(y))}{d_H(y)}. $$
Therefore, using Lemma \ref{l:3.5} for $B_r$ and the half plane $H:= \{\ell >0\}$ for each line, we conclude
\begin{align*}
{\rm {\rom{3}}}&=\frac{c}{r} \int_{{B_r} \cap \big( B_1(x_0) \backslash B_{\rho/2}(x_0) \big)} \frac{V(d_{B_r}(y))}{d_{B_r}(y)}\frac{J(1)}{|y-x_0|^{n-2} \varphi(|y-x_0|)}dy
\\ & \quad  +  \frac{c}{r}\int_{  H \cap \big(  B_1(x_0) \backslash B_{\rho/2}(x_0)  \big) } \frac{V(\ell(y))}{\ell(y)}\frac{J(1)}{|x-y|^{n-2} \varphi(|x_0-y|)}dy\\
&\le \frac{c_{23}}{r} \frac{r}{V(r)} + \frac{c_{24}}{r}\int_{  H \cap \big(  B_1(x_0) \backslash B_{\rho/2}(x_0)  \big) } \frac{V(d_H(y))}{d_H(y)}\frac{1}{|x-y|^{n-2} \varphi(|x_0-y|)}dy  \le \frac{c_{25}}{V(r)}.
\end{align*}
 Combining estimates of \rom{1},\rom{2} and \rom{3} we arrive
 $$ |L(V(\psi))(x_0)| \le {\rm {\rom{1}}}+{\rm {\rom{2}}}+{\rm {\rom{3}}} \le (c_{18}+c_{20}+c_{25}) \frac{1}{V(r)} $$
and \eqref{e:5.2} follows. \qed

\subsection{Subsolution and Harnack inequality}\label{s:Harnack}
In this section we construct a subsolution from the barrier we have obtained in Proposition \ref{p:barrier1}. Recall that we defined the domain of infinitesimal generator $A$ by
$$ \DD= \DD(D)= \{ u \in C_0(D) : Au \in C(D)   \} $$
in Section \ref{s:MH}. It is uncertain whether $V(\psi) \in \DD(D)$ since $A(V(\psi))$ is not continuous in general. To make our barrier included in the domain of operator, we construct a new domain of generator which contains $V(\psi)$. For given $C^{1,1}$ bounded open set $D$ and open subset $U$ in $D$, define
$$ \FF= \FF(D,U):= \{ u \in C_0(D) : Au \in L^\infty(U)  \}. $$
for the usage of proof. Denote $\FF(D)=\FF(D,D)$. Clearly $\FF(D,U_2) \subset \FF(D,U_1)$ for any $U_1 \subset U_2$. We first prove that $V(\psi) \in \FF(D)$.
\begin{lemma}\label{l:FF}
Let $\psi$ be the regularized version of $d_D$. Then, $A(V(\psi))= L(V(\psi))$ in $D$. Moreover, $V(\psi) \in \FF(D)$.
\end{lemma}
\pf Let $u \in C_0(D)$ be a twice-differentiable function in $D$. Assume that $\nabla^2 u$ is bounded in some $U \subset \subset D$. We first claim that
\begin{equation}\label{e:AL} Lu(x) = Au(x) \quad \mbox{for any} \quad x \in U. \end{equation}
Indeed, fix $x \in U$ and let $r_x>0$ be a constant satisfying $B=B(x,r_x) \subset U$. Without loss of generality we can assume $r_x \le 1$. Note that there exists a constant $c_1>0$ such that $2|u| + r_x^2\Ve \nabla^2 u\Ve \le c_1$ in $U$. Then we have
\begin{align}\begin{split}\label{e:AL1}
Au(x) &= \lim_{t \downarrow 0} \frac{P_tu(x) - u(x)}{t} = \lim_{t \downarrow 0} \frac{1}{t} \left(\int_{\R^n} u(x+y) p(t,|y|)dy - u(x)\right) \\
&= \lim_{t \downarrow 0} \int_{\R^n} \left( \frac{u(x+y) + u(x-y) }{2} -u(x) \right) \frac{p(t,|y|)}{t} dy.\end{split}
\end{align}
Since there is a constant $c_2>0$ such that $\frac{p(t,r)}{t} \le c_2 J(r)$ for any $t>0$ and $r>0$, we have
\begin{align*}
&\int_{\R^n} \left| \frac{u(x+y) + u(x-y) }{2} -u(x) \right| \frac{p(t,|y|)}{t} dy \\ &\le \int_{B} \left| \frac{u(x+y) + u(x-y) }{2} -u(x) \right| \frac{p(t,|y|)}{t} dy + \int_{B^c} \left| \frac{u(x+y) + u(x-y) }{2} -u(x) \right| \frac{p(t,|y|)}{t} dy \\ &\le c_1 \int_{B} \frac{|y|^2}{r_x^2} \frac{p(t,|y|)}{t}dy + c_1\int_{B^c}  \frac{p(t,|y|)}{t}dy \le c_1 \int_{\R^n} (\frac{|y|^2}{r_x^2} \land 1) (c_3J(|y|))dy< \infty
\end{align*}
for any $t>0$ so that we can apply dominate convergence theorem in the right-handed side of \eqref{e:AL1}.  Thus, using $\displaystyle \lim_{t \downarrow 0}$ $\frac{p(t,r)}{t} = J(r)$ we obtain
\begin{align*}Au(x) &= \lim_{t \downarrow 0} \int_{\R^n} \left( \frac{u(x+y) + u(x-y) }{2} -u(x) \right) \frac{p(t,|y|)}{t} dy \\ &= \int_{\R^n}\left( \frac{u(x+y) + u(x-y) }{2} -u(x) \right) J(|y|) dy = Lu(x). \end{align*}
This concludes the claim. Now, by Lemma \ref{l:V} we have that $V(\psi) \in C_0(D)$ is twice-differentiable and $\nabla^2 V(\psi)$ is locally bounded on $D$. Therefore, we arrive $L(V(\psi))= A(V(\psi))$ in $D$. It immediately follows from \eqref{e:5.1} that $V(\psi) \in \FF(D)$. \qed \\
Now we are ready to construct a subsolution with respect to the generator $A$.
\begin{lemma}[subsolution]\label{l:sub}
There exist a constant $C_4=C_4(n,a_1, a_2, \alpha_1, \alpha_2) > 0$ independent of $r$ and a radial function $w= w_r \in \FF(B_{4r})$ satisfying
		\begin{align*}
		\begin{cases}
		Aw \geq 0	&\text{in} ~ B_{4r} \setminus B_{r}, \\
		w \le V(r)	&\text{in} ~ B_{r}, \\
		w \geq C_4 V( 4r - \ve x \ve)	&\text{in} ~ B_{4r} \setminus B_{r}, \\
		w \equiv 0	&\text{in} ~ \Rn \setminus B_{4r},
		\end{cases}
		\end{align*}
		where $B_r := B(0, r)$.	
\end{lemma}
\pf Let $\Psi = \Psi_{4r}$ be the regularized version of $d_{B_{4r}}$ in \eqref{e:psi2} and choose a function $\eta \in C_c^\infty(B_1)$ satisfying $\Ve \eta \Ve_{C(B_1)}=1$ and $\eta \equiv 1$ on $B_{1/2}$. Define $\eta_r(x):=V(r) \eta(x/r) \in C_c^\infty(B_{r})$. Then, we have
\begin{align*}
|A\eta_r (x) |
&= |L\eta_r(x)| \le \int_{\Rn} \left| \frac{\eta_r(x+y)+\eta_r(x-y)}{2} - \eta_r(x) \right| J(|y|) dy \\
&\le \left( \Ve \nabla^2 \eta_r \Ve_{L^\infty(B_r)} + \Ve \eta_r \Ve_{L^\infty(B_r)} \right)  \int_{\Rn} \big( |y|^2 \land 1 \big) J(|y|)dy < \infty
\end{align*}
for any $x \in \R^n$, which implies $\eta_r \in \FF(B_{4r})$. Also, for $x \in B_{4r} \backslash B_{r}$,
\begin{align*}
A\eta_r (x) 
&= L\eta_r(x) = \int_{\Rn} \frac{\eta_r(x+y)+\eta_r(x-y)}{2} \frac{J(1)}{|y|^n\varphi(|y|)} dy \\
&= \int_{\Rn} \eta_r (x+y) \frac{J(1)}{|y|^n\varphi(|y|)} dy \ge \int_{ B(-x,r/2 ) }  \frac{V(r) J(1)}{|y|^n\varphi(|y|)} dy \ge \frac{ c_1(r/2)^n V(r)}{(9r/2)^n \varphi(9r/2)} \ge \frac{c_2}{V(4r)}.
\end{align*}
Here we used \eqref{e:V-asymp} and \eqref{e:V-wsc} for the last inequality. 

Define a function $\tilde{w}_r$ by
\begin{align*}
\tilde{w}_r = \frac{c_2}{C_3} V(\Psi) + \eta_r,
\end{align*}
where $C_3$ is the constant in Proposition \ref{p:barrier1}. We have $\tilde{w}_r \in \FF(B_{4r})$ by Lemma \ref{l:FF}. Also, for $x \in B_{4r} \backslash B_{r}$, using Proposition \ref{p:barrier1} and Lemma \ref{l:FF} again, we have
\begin{align*}
A\tilde{w}_r(x) = \frac{c_2}{C_3} AV(\Psi)(x) + A\eta_r(x) \ge - \frac{c_2}{C_3} |LV(\Psi)(x)| + A\eta_r(x) \ge -\frac{c_2}{V(4r)} + \frac{c_2}{V(4r)} = 0 
\end{align*}
and
\begin{align*}
\tilde{w}_r(x) = \frac{c_2}{C_3} V(\Psi(x)) \ge c_3 V(d_D(x)) = c_3 V(4r -|x|).
\end{align*}
For $x \in B_{r}$, 
$$ \tilde{w}_r(x) \le \frac{c_2}{C_3}V(4C r) + V(r) \le c_4 V(r) $$
by \eqref{e:psi2} and \eqref{e:V-wsc}.
Define $w_r(x) := \frac{1}{c_4} \tilde{w}_r(x)$. Then $w_r$ satisfies all assertions in Lemma \ref{l:sub} with constant $C_4 = \frac{c_3}{c_4}$, which is independent of $r$. \qed

We end this section with the Harnack inequality and the maximum principle of probabilistic version. For local operators, the Harnack inequality implies H\"older regularity of solutions of differential equations. However for nonlocal operators, as Silvestre mentioned in \cite{Sil}, this is not true because the nonnegativity of the function $u$ is required in the whole space $\Rn$. The Harnack inequality, maximum principle, and the subsolution constructed in Lemma \ref{l:sub} will play a key role in the proof of Theorem \ref{t:2}. We emphasize that the following theorem is the Harnack inequality for harmonic function with respect to $A$, and it does not imply the Harnack inequality for the viscosity solution with respect to $L$. See \cite{CS2} for the statement of Harnack inequality for viscosity solution.

\begin{theorem} [Harnack inequality] \cite[Theorem 2.2]{SV} \label{t:Harnack}
	Let $D$ be a bounded $C^{1,1}$ open set. Then, there exists a constant $C > 0$ such that for any ball $B(x_0,r) \subset D$, and any nonnegative function $u \in \FF(D)$ satisfying $Au = 0$ a.e. in $B(x_0, r)$, we have
	\begin{align*}
	\sup_{B(x_0, r/2)} u \leq C \inf_{B(x_0, r/2)} u.
	\end{align*}
\end{theorem}

Also, we have the following maximum principle. 

\begin{lemma}[Maximum principle]\label{l:m} Let $D$ be a bounded $C^{1,1}$ open set and $U$ be an open subset of $D$. If the function $u \in \FF(D,U)$ satisfies $Au=0$ a.e. in $U$ and $u \ge 0$ in $U^c$, then $u \ge 0$ in $\R^n$.
\end{lemma}
\pf Suppose that there exists $x \in U$ satisfying $u(x) < 0$. Since $u \in C_0(D)$, the set $U_- := \{ x \in \R^d : u(x) < 0 \}$ is bounded and open set with positive Lebesgue measure. For any $t>0$ we have
\begin{align*}
 \int_{U_-} P_t u(x)  - u(x)   dx &= \int_{U_-} \int_{\R^d} u(y) p(t,|x-y|) dy dx - \int_{U_-} u(x)dx \\
 &= \int_{\R^d} u(y) \int_{U_-} p(t,|x-y|)dx  dy - \int_{U_-} u(y)dy \\
 &= \int_{U_-^c} u(y)  \int_{U_-} p(t,|x-y|)dx  dy + \int_{U_-} u(y) \left( \int_{U_-} p(t,|x-y|)dx - 1 \right) dy \\
 &\ge \int_{U_-} u(y) \left( \int_{U_-} p(t,|x-y|)dx - 1 \right) dy.
\end{align*}
Since $U_-$ is bounded, $\rm{diam}(U_-)=:R<\infty$. Thus, for any $y \in U_- \subset B(y,R)$,
\begin{align*}
\frac{1-\int_{U_-} p(t,|x-y|)dx}{t} \ge \frac{1-\int_{B(y,R)} p(t,|x-y|)dx }{t} = \frac{1}{t} \left( 1 - \P^y (X_t \in B(y,R)) \right) = \frac{\P^0 (|X_t| \ge R)}{t}.
\end{align*}
Using heat kernel estimates in \cite[Theorem 21]{BGR1}, we have $p(t,r) \asymp  \left(\varphi^{-1}(t)^{-n} \land \frac{t}{r^n \varphi(r)}\right)$ for $(t,r) \in (0,1] \times \R_+$. Note that $\frac{t}{r^n\varphi(r)} \le \varphi^{-1}(t)^{-n}$ for $t  \le \varphi(r) $. Thus, there exists $\eps = \eps(R)>0$ satisfying 
$$ \frac{\P^0(|X_t| \ge R)}{t} \ge \frac{1}{t}\int_{R \le |z| \le 2R} p(t,|z|)dz \ge c \int_R^{2R} \frac{1}{r\varphi(r)}dr \ge  \eps \quad \mbox{for all} \quad t \in (0,\varphi(R)]. $$
Combining above estimates we obtain
$$ \int_{U_-} \frac{P_t u(x)- u(x)}{t} dx \ge -\eps \int_{U_-} u(y) dy \quad \mbox{for all} \quad t \in (0,\varphi(R)]. $$
Letting $t \to 0$, we conclude
$$ 0= \int_{U_-} Au(x)dx = \lim_{t \to 0}  \int_{U_-} \frac{P_t u(x)- u(x)}{t} dx \ge -\eps \int_{U_-} u(y)dy>0, $$
which is contradiction. Therefore, $u \ge 0$ in $\R^n$.
\qed

\subsection{Proof of Theorem \ref{t:2}} \label{s:pf of Thm 1.2}

In this section we will prove Theorem \ref{t:2}. More precisely, we prove the H\"older regularity for the function $u/V(d_D)$ up to the boundary of $D$. We will control the oscillation of this function using the Harnack inequality, the maximum principle and the subsolution constructed in Lemma \ref{l:sub}.

Let us adopt notations in \cite[Definition 3.3]{ROS1}. Let $\kappa > 0$ be a fixed small constant and let $\kappa' = 1/2 + 2\kappa$. Given $x_0 \in \partial D$ and $r>0$, define
\begin{align*}
D_r = D_r(x_0) = B(x_0, r) \cap D
\end{align*}
and
\begin{align*}
D_{\kappa'r}^+ = D_{\kappa'r}^+(x_0) = B(x_0, \kappa'r) \cap \lb x \in D : -x \cdot \nu(x_0) \geq 2\kappa r \rb,
\end{align*}
where $\nu(x_0)$ is the unit outward normal at $x_0$. Since $D$ is a bounded $C^{1,1}$ open set, there exists $\rho_0 > 0$ such that for each $x_0 \in \partial D$ and $r \leq \rho_0$, there exists an orthonormal system $CS_{x_0}$ with its origin at $x_0$ and a $C^{1,1}$-function $\Psi : \R^{n-1} \rightarrow \R$ satisfying $\Psi(\tilde{0}) = 0, \nabla_{CS_{x_0}} \Psi(\tilde{0}) = 0, \Ve \Psi \Ve_{C^{1,1}} \leq \kappa$, and
\begin{align*}
\lb y=(\tilde{y}, y_n) ~ \text{in} ~ CS_{x_0} : \ve \tilde{y} \ve < 2r, \Psi(\tilde{y}) < y_n < 2r \rb \subset D.
\end{align*}
Then we have
\begin{align} \label{e:Dr1}
B(y, \kappa r) \subset D_r(x_0) ~~ \text{for all} ~~ y \in D^+_{\kappa' r}(x_0),
\end{align}
and we can take a $C^{1.1}$ subdomain $D_r^{1,1}$ satisfying $D_r \subset D_r^{1,1} \subset D_{2r}$ and
\begin{align} \label{e:Dr11}
\dist (y, \partial D_r^{1,1}) = d_D(y)
\end{align}
for all $y\in D_r$. Since $D_r$ is not $C^{1,1}$ in general, we will use this subdomain instead of $D_r$.

Since $D$ is bounded and $C^{1,1}$ again, we can assume that for each $x_0 \in \partial D$ and $r \leq \rho_0$,
\begin{align} \label{e:Dr2}
B( y^* - 4\kappa r \nu(y^*), 4\kappa r) \subset D_r(x_0) ~~\text{and} ~~ B(y^* - 4\kappa r \nu(y^*), \kappa r) \subset D_{\kappa' r}^+ (x_0)
\end{align}
for all $y \in D_{r/2}(x_0)$, where $y^* \in \partial D$ is the unique boundary point satisfying $\ve y - y^* \ve = d_D(y)$. 

The following oscillation lemma is the key lemma to prove Theorem \ref{t:2}.

\begin{lemma} [Oscillation lemma] \label{l:osc}
Assume $f \in C(D)$ and let $u \in \DD$ be the viscosity solution of \eqref{e:pde1}. Then there exist constants $\gamma \in (0, 1)$ and $C_1 > 0$, depending only on $n, a_1, a_2, \alpha_1, \alpha_2$ and $D$, such that
\begin{align} \label{e:osc}
\sup_{D_r(x_0)} \frac{u}{V(d_D)} - \inf_{D_r(x_0)} \frac{u}{V(d_D)} \leq C_1 V(r)^\gamma \Ve f \Ve_{L^\infty(D)} 
\end{align}
for any $x_0 \in \partial D$ and $r > 0$.
\end{lemma}

To prove the oscillation lemma, we need some preparation. Note that in the following two lemmas we aim to verify inequalities for every function $u \in \FF$, since we want to utilize the subsolution constructed in Lemma \ref{l:sub}. The first one is a generalized version of Harnack inequality.

\begin{lemma} [Harnack inequality] \label{l:harnack}
There exists a constant $C_2 = C_2(n, a_1, a_2, \alpha_1, \alpha_2, D) > 0$ such that for any $r \leq \rho_0, x_0 \in \partial D$ and nonnegative function $u \in \FF(D,D_r^{1,1})$,
\begin{align} \label{e:Harnack}
\sup_{D^+_{\kappa' r}(x_0)} \frac{u}{V(d_D)} \leq C_2 \left( \inf_{D^+_{\kappa' r}(x_0)} \frac{u}{V(d_D)} + \Ve Au \Ve_{L^\infty(D_r^{1,1})} V(r) \right).
\end{align}
\end{lemma}

\pf
We first prove that if a nonnegative function $v$ satisfies $Av = 0$ a.e. in $D^{1,1}_r$, then 
\begin{align} \label{e:h for v/V}
\sup_{D^+_{\kappa' r}(x_0)} \frac{v}{V(d_D)} \leq c \inf_{D^+_{\kappa' r}(x_0)} \frac{v}{V(d_D)}
\end{align}
for a constant $c>0$ which is independent of $r$ and $v$. Indeed, for each $y \in D^+_{\kappa' r}$, we have $B(y, \kappa r) \subset D^{1,1}_r$ by \eqref{e:Dr1} hence $Av = 0$ a.e. in $B(y, \kappa r)$. We may cover $D^+_{\kappa' r}$ by finitely many balls $B(y_i, \kappa r/2)$. Here the number of balls is independent of $r$. By the Theorem \ref{t:Harnack}, we have for each $i$,
\begin{align*}
\sup_{B(y_i, \kappa r/2)} v \leq c_1 \inf_{B(y_i, \kappa r/2)} v.
\end{align*}
If $x \in B(y_i, \kappa r/2)$, we have $\kappa r /2 \leq d_D(x) \leq r/2 + 5\kappa r / 2$. Thus, using \eqref{e:V-wsc} we obtain
\begin{align*}
\sup_{B(y_i, \kappa r/2)} \frac{v}{V(d_D)} \leq \sup_{B(y_i, \kappa r/2)} \frac{v}{V(\kappa r / 2)} \leq c_2 \inf_{B(y_i, \kappa r/2)} \frac{v}{V(r/2 + 5\kappa r/2)} \leq c_2 \inf_{B(y_i, \kappa r/2)} \frac{v}{V(d_D)}.
\end{align*}
Now \eqref{e:h for v/V} follows from the standard covering argument, possibly with a larger constant.

We next prove \eqref{e:Harnack}. Let us write $u = u_1 + u_2$, where $u_1 := u + R^{D_r^{1,1}}Au$ and $u_2 := -R^{D_r^{1,1}} Au$. We claim that $u_1 \ge 0$ in $\Rn$ and $Au_1 = 0$ a.e. in $D^{1,1}_r$.

Following the calculations of \eqref{e:A} we obtain that for any open subset $U \subset D$, $x \in U$ and $u \in \FF(D,U)$,
\begin{equation}\label{e:A1}
Au(x) = \lim_{ t \downarrow 0} \frac{P_t u(x) - u(x)}{t}=\lim_{ t \downarrow 0} \frac{P_t^U u(x) - u(x)}{t}. 
\end{equation}
Let us emphasize that we only have used $u \in C_0(D)$ in \eqref{e:A} so we can repeat the same argument for $u \in \FF(D,U)$.
%\red{------------------------------- This is draft of proof---------------------------------------------}
%
%We first show $u_1 \ge 0$. Let $-A = \int_0^\infty \lambda dE_\lambda$. Then we have $P^U_t = \int_0^\infty e^{-\lambda t} dE_\lambda$. Since $R^U$ is positive operator we have 
%\begin{align*}
%R^U Au(x) \ge R^U A^U u(x) = A^U R^U u(x) = u(x)
%\end{align*}
%where we used \eqref{e:Au} and $u \in \FF$ for the last equality. For the second equality we observe that
%\begin{align*}
% R^U A^U u(x) &= \lim_{ N \to \infty} \int_0^N P_t^U A^U u(x) dt = \lim_{N \to \infty} \int_0^N A^U P_t^U u(x) dt = \lim_{N \to \infty} A^U \left(\int_0^N P_t^U u(x)dt \right) \\
% &= 
%\end{align*}
%

Let $g \in L^\infty(U)$. Deducing $R^U g \in C_0(U)$ from Proposition \ref{p:R} and \eqref{d:R}, we obtain the following counterpart of \eqref{e:Au}: For any $x \in U$,
\begin{align}
\begin{split}\label{e:A2}
AR^U g(x) &= A \left(\int_0^\infty P_s^U g(\cdot )      ds \right)(x) = \lim_{t \downarrow 0} \frac{1}{t} \left( P^U_t\left(\int_0^\infty P_s^U g(\cdot )      ds \right)(x) - \int_0^\infty P_s^U g(x)      ds  \right) \\
& = \lim_{t \downarrow 0} \frac{1}{t} \left( \int_0^\infty P_{s+t}^U g(x) ds - \int_0^\infty P_s^U g(x)ds \right) \\ &= -\lim_{t \downarrow 0} \frac{ \int_0^t P_s^U g(x)ds}{t} =-\lim_{t \downarrow 0} \frac{ \int_0^t P_s g(x)ds}{t}.
\end{split}
\end{align}
Here we used \eqref{e:A1} for the first line. Let 
$$U_g := \{ x \in U : \lim_{r \downarrow 0} \frac{1}{r^n} \int_{B(x,r)} |g(x)-g(y)|dy =0  \}. $$
Then, we have $| U \setminus U_g |=0$ since $g \in L^\infty(U) \subset L^1(U)$. For $x \in U_g$, we have
$$|P_t g(x) - g(x) | = \lv \int_{\R^n} p(t,|x-y|)(g(y)-g(x))dy \rv \le \int_{\R^n} p(t,|x-y|) |g(y)-g(x)|dy. $$
Let $\eps>0$. Using $p(t,r) \asymp \left( \varphi^{-1}(t)^{-n} \land \frac{t}{r^n \varphi(r)}\right)$ for $t \in (0,1] \times \R_+$ in \cite[Theorem 21]{BGR1} again, there exist constants $c_3(\eps), c_4(\eps)>0$  such that for any $t \in (0,1]$ and $r>0$,
$$ p(t,r) \le c_4 \varphi^{-1}(t)^{-n} $$
and 
$$ \P^x (|X_t| > c_3 \varphi^{-1}(t)) \le \eps. $$
Indeed, using \eqref{e:varphi-inf} and \eqref{e:varphi-wsc} we have
$$\P^x (|X_t| > c_3 \varphi^{-1}(t)) = \int_{|z|>c_3 \varphi^{-1}(t)} p(t,|z|)dz \le c_4 t\int_{c_3 \varphi^{-1}(t)}^\infty \frac{dr}{r\varphi(r)} \le \frac{c_5t} {\varphi(c_3\varphi^{-1}(t))} \le c_6c_3^{-2\alpha_1}. $$
Thus, we obtain
\begin{align*}
|P_t g(x) - g(x) | &\le \int_{B(x, c_3\varphi^{-1}(t))} p(t,|x-y|) |g(y)-g(x)| dy + \int_{B(x, c_3\varphi^{-1}(t))^c} p(t,|x-y|) |g(y)-g(x)|dy \\
&\le c_4 \varphi^{-1}(t)^{-n} \int_{B(x, c_3 \varphi^{-1}(t))} |g(y)-g(x)|dy + 2 \Ve g \Ve_{\infty} \int_{B(x, c_3\varphi^{-1}(t))} p(t,|x-y|)dy \\
&\le c_4 \varphi^{-1}(t)^{-n} \int_{B(x, c_3 \varphi^{-1}(t))} |g(y)-g(x)|dy + 2\Ve g \Ve_{\infty} \eps.
\end{align*}
Since $\eps>0$ is arbitrary and $x \in U_g$, we conclude
$$ \lim_{t \downarrow 0} |P_t g(x)- g(x)| =0. $$
Combining this with \eqref{e:A2} we arrive that for any open subset $U \subset D$ and $g \in L^\infty(D)$,
\begin{equation}\label{e:g}
 AR^{U} g= -g \quad \mbox{a.e. in} \quad U. 
\end{equation}
Since $u \in \FF(D,U)$, we have $Au \in L^\infty(U)$. Thus, taking $U=D_r^{1,1}$ and $g=Au$ in \eqref{e:g} we conclude
$$ Au_1 = Au + AR^{D_r^{1,1}}Au = 0 \quad \mbox{a.e. in } D_r^{1,1}. $$
Also, $u_1 \ge 0$ follows from applying Lemma \ref{l:m} with above equation and $u_1 = u \ge 0$ in $\R^n \backslash D_r^{1,1}$.

Applying \eqref{e:h for v/V} to $u_1$, we get
\begin{align*}
\sup_{D^+_{\kappa' r}} \frac{u_1}{V(d_D)} \leq c_7 \inf_{D^+_{\kappa' r}} \frac{u_1}{V(d_D)}.
\end{align*}
Meanwhile, using \eqref{e:Dr11} and Lemma \ref{l:u} we have
\begin{align*}
\ve u_2(x) \ve \leq c_8 \Ve Au \Ve_{L^\infty(D_r^{1,1})} V(\diam(D_r^{1,1})) V(\dist(x, \partial D^{1,1}_r)) \leq c_9 \Ve Au \Ve_{L^\infty(D_r^{1,1})} V(r) V(d_D(x))
\end{align*}
for all $x \in D_r^{1,1}$. Therefore, combining above two inequalities we conclude that
\begin{align*}
\sup_{D^+_{\kappa' r}} \frac{u}{V(d_D)} 
&\leq \sup_{D^+_{\kappa' r}} \frac{u_1}{V(d_D)} + \sup_{D^+_{\kappa' r}} \frac{u_2}{V(d_D)} \leq c_5 \inf_{D^+_{\kappa' r}} \frac{u_1}{V(d_D)} + \sup_{D^+_{\kappa' r}} \frac{u_2}{V(d_D)} \\
&\leq c_5 \inf_{D^+_{\kappa' r}} \frac{u}{V(d_D)} + (c_5 + 1) \sup_{D^+_{\kappa' r}} \frac{\ve u_2 \ve}{V(d_D)} \leq C_2 \left( \inf_{D^+_{\kappa' r}} \frac{u}{V(d_D)} + \Ve Au \Ve_{L^\infty(D_r^{1,1})} V(r) \right).
\end{align*}
 \qed

The next lemma gives the link between $D^+_{\kappa' r}$ and $D_{r/2}$. Here we are going to use the subsolution $w$ in Lemma \ref{l:sub}.

\begin{lemma} \label{l:4.9}
Let $r \leq \rho_0, x_0 \in \partial D$. If $u \in \FF(D,D_r^{1,1})$ is nonnegative, then there exists a constant $C_3 = C_3(n, a_1, a_2, \alpha_1, \alpha_2, D) > 0$ such that
\begin{align*}
\inf_{D^+_{\kappa' r}(x_0)} \frac{u}{V(d_D)} \leq C_3 \left( \inf_{D_{r/2}(x_0)} \frac{u}{V(d_D)} + \Ve Au \Ve_{L^\infty(D_r^{1,1})} V(r) \right).
\end{align*}
\end{lemma}

\pf
First assume that $Au$ is nonnegative. As in the proof of Lemma \ref{l:harnack}, we write $u = u_1 + u_2$, where $u_1 = u + R^{D_r^{1,1}} Au$ and $u_2 = -R^{D_r^{1,1}} Au$. Then $u_1$ is a nonnegative solution for
\begin{align*}
\begin{cases}
Au_1 = 0	&\text{a.e. in} ~ D^{1,1}_r, \\
u_1 = u	&\text{in} ~  \Rn \setminus D^{1,1}_r.
\end{cases}
\end{align*}
Let
\begin{align*}
m:= \inf_{D^+_{\kappa' r}} \frac{u_1}{V(d_D)} \geq 0.
\end{align*}
For $y \in D_{r/2}$, we have either $y\in D^+_{\kappa' r}$ or $d_D(y) < 4\kappa r$ by \eqref{e:Dr2}. 

If $y\in D^+_{\kappa' r}$, then clearly
\begin{align} \label{e:inf}
m \leq \frac{u_1(y)}{V(d_D(y))}.
\end{align}

If $d_D(y) < 4 \kappa r$, let $y^*$ be the closest point to $y$ on $\partial D_r^{1,1}$ and let $\tilde{y} = y^* - 4\kappa r \nu (y^*)$. By \eqref{e:Dr2}, we have $B_{4\kappa r}(\tilde{y}) \subset D_r$ and $B_{\kappa r}(\tilde{y}) \subset D^+_{\kappa' r}$. 

Now consider $w \in \FF(B_{4\kappa r} (\tilde{y})) \subset \FF(D,B_{4\kappa r} (\tilde{y}) \backslash B_{\kappa r}(\tilde{y}))$ satisfying
\begin{align*}
\begin{cases}
Aw \geq 0	&\text{in} ~ B_{4\kappa r} (\tilde{y}) \setminus B_{\kappa r}(\tilde{y}), \\
w \leq V(\kappa r) &\text{in} ~ B_{\kappa r} (\tilde{y}), \\
w \geq c_1 V(4\kappa r - \ve x - \tilde{y} \ve) &\text{in} ~ B_{4\kappa r} (\tilde{y}) \setminus B_{\kappa r} (\tilde{y}), \\
w \equiv 0	&\text{in} ~ \Rn \setminus B_{4\kappa r} (\tilde{y}),
\end{cases}
\end{align*}
which can be obtained by translating the subsolution in Lemma \ref{l:sub}. Since $Au_1 = 0$ a.e. in $B_{4\kappa r} (\tilde{y})$,  we have
\begin{align*}
\begin{cases}
Au_1 = 0 \leq A(mw) &\text{a.e. in} ~ B_{4\kappa r} (\tilde{y}) \setminus B_{\kappa r}(\tilde{y}), \\
u_1 \geq mV(d_D) \geq mw &\text{in} ~ B_{\kappa r} (\tilde{y}), \\
u_1 \geq 0 = m w &\text{in} ~ \Rn \setminus B_{4\kappa r} (\tilde{y}).
\end{cases}
\end{align*}
Now by the maximum principle in Lemma \ref{l:m} with the function $u_1 -mw$ and $U=B_{4\kappa r} (\tilde{y}) \setminus B_{\kappa r}(\tilde{y})$, we obtain $u_1 \geq mw$ in $\Rn$. In particular, for $y \in B_{4\kappa r} (\tilde{y}) \setminus B_{\kappa r}(\tilde{y})$,
\begin{align*}
u_1(y) \geq c_1 m V(4\kappa r - \ve y - \tilde{y} \ve ) = c_1 m V(d_D(y)).
\end{align*}
Therefore, we obtain
\begin{align*}
\inf_{D^+_{\kappa' r}} \frac{u_1}{V(d_D)} \leq c_2 \inf_{D_{r/2}} \frac{u_1}{V(d_D)}.
\end{align*}
On the other hand, $u_2$ satisfies
\begin{align*}
\ve u_2(x) \ve \leq c_3 \Ve Au \Ve_{L^\infty(D_r^{1,1})} V(r) V(d_D(x))
\end{align*}
for all $x \in D_r^{1,1}$, which gives the desired result. \qed

We prove the oscillation lemma \eqref{e:osc} by using Lemmas \ref{l:harnack} and \ref{l:4.9}.
\\

\noindent{\bf Proof of Lemma \ref{l:osc}}
As a consequence of Remark \ref{r:3.7}, by dividing $\Ve f \Ve_{L^\infty(D)}$ on both sides of \eqref{e:pde1} if necessary, we may assume $\Ve f \Ve_{L^\infty (D)} \leq 1$ and $\Ve u \Ve_{C(D)}=\Ve R^Df \Ve_{C(D)} \leq c_1$ without loss of generality. Fix $x_0 \in \partial D$. We will prove that there exist constants $c_2 > 0, \rho_1 \in (0, \rho_0 /16]$, and $\gamma \in (0,1)$
%, depending only on $n, a_1, a_2, \alpha_1, \alpha_2$, and $D$, 
and monotone sequences $(m_k)_{k\geq 0}$ and $(M_k)_{k \geq 0}$ such that $M_k - m_k = V(r_{k+1} / 2)^\gamma,$
\begin{align*}
-V(\rho_1 / 16) \leq m_k \leq m_{k+1} < M_{k+1} \leq M_k \leq V(\rho_1 / 16),
\end{align*}
and
\begin{align*}
m_k \leq \frac{u}{c_2 V(d_D)} \leq M_k ~~ \text{in} ~ D_{r_k} = D_{r_k}(x_0)
\end{align*}
for all $k \geq 0$, where $r_k = \rho_1 8^{-k}$. If we have such constants and sequences, then for any $0< r \leq \rho_1$ we have $k \ge 0$ satisfying $r \in (r_{k+1},r_k]$ and
\begin{align*}
\sup_{D_r} \frac{u}{V(d_D)} - \inf_{D_r} \frac{u}{V(d_D)} \leq \sup_{D_{r_k}} \frac{u}{V(d_D)} - \inf_{D_{r_k}} \frac{u}{V(d_D)} \leq c_2 (M_k-m_k) = c_2 V(r_{k+1} / 2)^\gamma \leq c_2 V(r)^\gamma.
\end{align*}
Also, for any $r > \rho_1$ we have 
\begin{align*}
\sup_{D_r} \frac{u}{V(d_D)} - \inf_{D_r} \frac{u}{V(d_D)} \leq c_3 \leq c_4 V(\rho_1)^\gamma \leq c_4 V(r)^\gamma
\end{align*}
by Lemma \ref{l:u}. Above two inequalities conclude the lemma so it suffices to construct such constants and sequences.

Let us use the induction on $k$. The case $k=0$ follows from Lemma \ref{l:u} provided we take $c_2$ large enough. The constants $\rho_1$ and $\gamma$ will be chosen later. Assume that we have sequences up to $m_k$ and $M_k$. Let $\psi$ be the regularized version of $d_D$. We may assume that $\psi = d_D$ in $\lb d_D(x) \leq \rho_1 \rb$. Define
\begin{align*}
u_k = V(\psi) \left( \frac{u}{c_2 V(\psi)} - m_k \right) = \frac{1}{c_2} u - m_k V(\psi)
\end{align*}
in $\Rn$. Note that $u_k \in \FF(D)$ since $Au=f$ by the consequence of Theorem \ref{t:sol}. Moreover, for $x \in D^{1,1}_{r_k / 4}$ we have $u_k^- \in C^2(x)$ since we know that $u_k^- \equiv 0$ in $B(x_0, r_k)$ by the induction hypothesis. Thus, we have $Au_k^-(x) = Lu_k^-(x)$ by \eqref{e:gen2}, which implies that $Au_k^-$ is well-defined in $D^{1,1}_{r_k / 4}$, and so is $Au_k^+$. We will apply Lemmas \ref{l:harnack} and \ref{l:4.9} for the function $u_k^+$ and $r = r_k / 4$ to find $m_{k+1}$ and $M_{k+1}$. By \eqref{e:5.1} and Lemma \ref{l:FF}, we have 
\begin{align} \label{e:Lu+}
\begin{split}
\ve Au_k^+ \ve 
&\leq \ve Au_k \ve + \ve Au_k^- \ve \leq \lv \frac{1}{c_2} Au - m_k AV(\psi) \rv + \ve Au_k^- \ve \\
&\leq \left(\frac{1}{c_2} \ve f \ve + V(\rho_1/16) \ve L(V(\psi)) \ve\right) + \ve Au_k^- \ve \leq c_3 + \ve Au_k^- \ve 
\end{split}
\end{align}
in $D$. Thus, we need to estimate $\ve Au_k^- \ve$ in $D_{r_k / 4}^{1,1}$ for the usage of Lemmas \ref{l:harnack} and \ref{l:4.9}.

Let $x \in D^{1,1}_{r_k / 4}$. By the induction hypothesis, we have $u_k^- \equiv 0$ in $B(x_0, r_k)$, which implies that $u_k^- \in C^2(x)$. Thus, we compute the value $Au_k^-(x)$ using the operator $L$ as follows:
\begin{equation} \label{e:Lu_k^-}
\begin{split}
0 \leq Au_k^-(x) = Lu_k^-(x) &= \frac{1}{2} \int_{\Rn} \left( u_k^-(x+h) + u_k^-(x-h) \right) \frac{J(1)}{\ve h \ve^n \varphi(\ve h \ve)} \, dh \\
&= \int_{x+h \notin B_{r_k}} u_k^-(x+h) \frac{J(1)}{\ve h \ve^n \varphi(\ve h \ve)} \, dh.
\end{split}
\end{equation}
For any $y \in B_{r_0} \setminus B_{r_k}$, there is $0 \le j<k$ such that $y\in B_{r_j} \setminus B_{r_{j+1}}$. Since $c_2^{-1} u \geq m_j V(\psi)$ and $d_D = \psi$ in $B_{r_j}$, we have
\begin{align*}
u_k(y) &= c_2^{-1} u(y) - m_k V(\psi(y)) \geq (m_j - m_k)V(\psi(y)) \\
&\geq (m_j - M_j + M_k - m_k) V(d_D(y)) \geq - (V(r_{j+1} / 2)^\gamma - V(r_{k+1} / 2)^\gamma ) V(r_j).
\end{align*}
It follows from $r_{j+1} \leq \ve y - x_0 \ve < r_j \leq 8\ve y - x_0 \ve \leq 1$ that
\begin{equation} \label{e:u_k^-}
\begin{split}
u_k^- (y) 
&\leq c_4 \left( V(\ve y - x_0 \ve / 2)^\gamma - V(r_k / 16)^\gamma \right) V(8 \ve y - x_0 \ve) \\
&\leq c_5 \left( V(\ve y - x_0 \ve / 2)^\gamma - V(r_k / 16)^\gamma \right) V( \ve y - x_0 \ve / 2).
\end{split}
\end{equation}
Note that \eqref{e:u_k^-} possibly with a larger constant also holds for $y \in \Rn \setminus B_{r_0}$ because $\Ve u_k \Ve_{C(\Rn)} \leq c_1 c_2^{-1} + V(1/16) V(\tilde{C})$ for any $k$ and
\begin{align*}
\left( V(\ve y - x_0 \ve / 2)^\gamma - V(r_k / 16)^\gamma \right) V(\ve y - x_0 \ve / 2) \geq \left( V(\rho_1 / 2)^\gamma - V(\rho_1 /16)^\gamma \right) V(\rho_1 / 2) > 0
\end{align*}
for any $y \in \Rn \setminus B_{r_0}$. Thus, by \eqref{e:Lu_k^-} and \eqref{e:u_k^-}, we have
\begin{align*}
\ve Au_k^-(x) \ve 
&\leq c_6 \int_{x+y \notin B_{r_k}} \left( V(\ve x+h - x_0 \ve / 2)^\gamma - V(r_k / 16)^\gamma \right) \frac{V(\ve x+y-x_0 \ve / 2)}{\ve h \ve^n \varphi(\ve h \ve)} \, dh.
\end{align*}
If $x+y \notin B_{r_k}$, then $\ve h \ve \geq \ve x+h-x_0 \ve - \ve x-x_0 \ve \geq r_k - r_k / 2 = r_k / 2$ and $\ve x+h-x_0 \ve \leq r_k / 2 + \ve h \ve \leq 2\ve h \ve$. Thus, recalling that $\sP_1(r)=\int_r^\infty \frac{ds}{s\varphi(s)}$, we obtain
\begin{align*}
\ve Au_k^-(x) \ve 
\leq c_6& \int_{\ve h \ve \geq r_k / 2} \left( V(\ve h \ve)^\gamma - V(r_k / 16)^\gamma \right) \frac{V(\ve h \ve)}{\ve h \ve^n \varphi(\ve h \ve)} \, dh \\
\leq c_7& \int_{r_k/2}^\infty \left( V(s)^\gamma - V(r_k / 16)^\gamma \right) V(s) d(-\mathcal{P}_1)(s) \\
= c_7& \bigg( \left[ - \left( V(s)^\gamma - V(r_k / 16)^\gamma \right) V(s) \mathcal{P}_1(s) \right]_{r_k/2}^\infty \\
&+ \int_{r_k / 2}^\infty \left( (1+\gamma)V(s)^\gamma - V(r_k / 16)^\gamma \right) V'(s) \mathcal{P}_1(s) ds \bigg) =: c_7 \left( {\rm \rom{1}} + {\rm \rom{2}} \right).
\end{align*}
By \eqref{e2} we have
\begin{align*}
\lim_{s\rightarrow \infty} \left( V(s)^\gamma -V(r_k / 16)^\gamma \right) V(s) \mathcal{P}_1(s) \leq c_8 \lim_{s\rightarrow\infty} \frac{V(s)^\gamma - V(r_k / 16)^\gamma}{V(s)} = 0,
\end{align*}
hence
\begin{align*} 
{\rm \rom{1}} \leq c_8 \frac{ V(r_k/2)^\gamma - V(r_k / 16)^\gamma }{V(r_k / 2)}.
\end{align*}
Also, using \eqref{e2} again we have
\begin{equation*} 
\begin{split}
{\rm \rom{2}} &\leq c_8 \int_{r_k / 2}^\infty \left( (1+\gamma)V(s)^\gamma - V(r_k / 16)^\gamma \right) \frac{V'(s)}{V(s)^2} \, ds \\
&= c_8 \left( \frac{1+\gamma}{1-\gamma} V(r_k/2)^\gamma - V(r_k/16)^\gamma \right) \frac{1}{V(r_k/2)}.
\end{split}
\end{equation*}
Therefore, combining above two inequalities and using \eqref{e:V-wsc} we get
\begin{align*}
\ve Au_k^-(x) \ve 
&\leq c_9 \left( \frac{2}{1-\gamma} V(r_k / 2)^\gamma - 2V(r_k / 16)^\gamma \right) \frac{1}{V(r_k / 2)} \\
&\leq c_9 \left( \frac{2}{1-\gamma} \left( c_{10} 64^{\alpha_2} \right)^\gamma - 2 \left( c_{10}^{-1} 8^{\alpha_1} \right)^\gamma \right) \frac{V(r_{k+2}/2)^\gamma }{V(r_k / 4)} \\
&=: c_9 \eps_\gamma \frac{V(r_{k+2}/2)^\gamma }{V(r_k / 4)}
\end{align*}
and hence
\begin{align*}
\Ve Au_k^+ \Ve_{L^\infty(D_{r_k / 4}^{1,1})} \leq c_{11} \left( 1+ \eps_\gamma \frac{V(r_{k+2}/2)^\gamma}{V(r_k / 4)} \right).
\end{align*}
Note that $\eps_\gamma \rightarrow 0$ as $\gamma \rightarrow 0$.

Now we apply Lemma \ref{l:harnack} and \ref{l:4.9} for $u_k^+ \in \FF(D,D_{r_k/4}^{1,1})$. Since $u_k = u_k^+$ and $d_D = \psi$ in $D_{r_k}$, we have
\begin{align*}
\sup_{D^+_{\kappa' r_k / 4}} \left( \frac{u}{c_2 V(\psi)} - m_k \right) 
&\leq c_{12} \left( \inf_{D^+_{\kappa' r_k / 4}} \left( \frac{u}{c_2 V(\psi)} - m_k \right) + V(r_k / 4) + \eps_\gamma V(r_{k+2} / 2)^\gamma \right) \\
&\leq c_{13} \left( \inf_{D_{r_{k+1}}} \left( \frac{u}{c_2 V(\psi)} - m_k \right) + V(r_k / 4) + \eps_\gamma V(r_{k+2}/2)^\gamma \right).
\end{align*}
Repeating this procedure with the function $u_k = M_k V(d_D) - c_2^{-1}u$ instead of $u_k = c_2^{-1} u - m_k V(d_D)$, we also have
\begin{align*}
\sup_{D^+_{\kappa' r_k / 4}} \left( M_k - \frac{u}{c_2 V(\psi)} \right) \leq c_{14} \left( \inf_{D_{r_{k+1}}} \left( M_k - \frac{u}{c_2 V(\psi)} \right) + V(r_k / 4) + \eps_\gamma V(r_{k+2}/2)^\gamma \right).
\end{align*}
Adding up these two inequalities, we obtain
\begin{align*}
M_k - m_k \leq c_{15} \left( \inf_{D_{r_{k+1}}} \frac{u}{c_2 V(\psi)} - \sup_{D_{r_{k+1}}} \frac{u}{c_2 V(\psi)} + M_k - m_k + V(r_k / 4) + \eps_\gamma V(r_{k+2} / 2)^\gamma \right).
\end{align*}
Thus, recalling that $M_k - m_k = V(r_{k+1} / 2)^\gamma$, we get
\begin{align*}
\sup_{D_{r_{k+1}}} \frac{u}{c_2 V(\psi)} - \inf_{D_{r_{k+1}}} \frac{u}{c_2 V(\psi)} 
&\leq \frac{c_{15}-1}{c_{15}} V(r_{k+1}/2)^\gamma + V(r_k / 4) + \eps_\gamma V(r_{k+2}/2)^\gamma \\
&\leq \left( \frac{c_{15}-1}{c_{15}} c_{16}^\gamma + c_{17}^\gamma V(\rho_1)^{1-\gamma} + \eps_\gamma \right) V(r_{k+2}/2)^\gamma.
\end{align*}
Now we choose $\gamma$ and $\rho_1$ small enough so that
\begin{align*}
\frac{c_{15}-1}{c_{15}} c_{16}^\gamma + c_{17}^\gamma V(\rho_1)^{1-\gamma} + \eps_\gamma \leq 1,
\end{align*}
and it yields that
\begin{align*}
\sup_{D_{r_{k+1}}} \frac{u}{c_2V(\psi)} - \inf_{D_{r_{k+1}}} \frac{u}{c_2 V(\psi)} \leq V(r_{k+2}/2)^\gamma.
\end{align*}
Therefore, we are able to choose $m_{k+1}$ and $M_{k+1}$.
\qed

Finally, we prove the Theorem \ref{t:2} using the Lemma \ref{l:osc}. \\

\noindent{\bf Proof of Theorem \ref{t:2}}
By Remark \ref{r:3.7}, by dividing $\Ve f \Ve_{L^\infty(D)}$ on both sides of \eqref{e:pde1} if necessary, we may assume that $\Ve f \Ve_{L^\infty(D)} \leq 1$ and $\Ve u \Ve_{C(D)} \leq c_1$. We first show that the following holds for any $x \in D$:
\begin{align*}
\left[ \frac{u}{V(d_D)} \right]_{C^\beta (B(x, r/2))} \leq \frac{C}{r^{\beta}V(r)}
\end{align*}
for each $0<\beta \leq \alpha_1$, where $r = d_D(x)$. We are going to use the inequality
\begin{align} \label{e:Cbeta}
\left[ \frac{u}{V(d_D)} \right]_{C^\beta} \leq \Ve u \Ve_C \left[ \frac{1}{V(d_D)} \right]_{C^\beta} + [u]_{C^\beta} \lV \frac{1}{V(d_D)} \rV_C.
\end{align}

From \eqref{e:CV} we know that $[u]_{C^V ( B(x, r/2))} \leq c_2$. Thus, we have $[u]_{C^\beta (B(x, r/2))} \leq c_3$ for each $0<\beta \leq \alpha_1$. Since $d_D(y) \geq r/2$ for $y \in B(x, r/2)$, we have
\begin{align*}
\lV \frac{1}{V(d_D)} \rV_{C (B(x, r/2))} \leq \frac{c_4}{V(r)}
\end{align*}
and
\begin{align*}
\left[ \frac{1}{V(d_D)} \right]_{C^{0,1}(B(x, r/2))} &\leq \sup_{y,z\in B(x, r/2)} \frac{\ve V(d_D(y))^{-1} - V(d_D(z))^{-1} \ve}{\ve y - z \ve} \\
&\leq \sup_{y, z \in B(x, r/2)} \frac{V'(d^*)}{V(d^*)^2} \frac{ \ve d_D(y) - d_D(z) \ve}{\ve y - z \ve} \\
&\leq c_5 \left( \sup_{y, z \in B(x, r/2)} \frac{1}{d^*V(d^*)} \right) [d]_{C^{0,1}(B(x, r/2))} \\
&\leq \frac{c_6}{rV(r)},
\end{align*}
where $d^*$ is a value in $[d_D(y), d_D(z)]$, so $d^* \geq r/2$. Thus, by interpolation, we obtain
\begin{align*}
\left[ \frac{1}{V(d_D)} \right]_{C^\beta (B(x, r/2))} &\leq c_7 \lV \frac{1}{V(d_D)} \rV_{C (B(x, r/2))}^{1-\beta} \left[ \frac{1}{V(d_D)} \right]_{C^{0,1}(B(x, r/2))}^{\beta} \leq \frac{c_8}{r^\beta V(r)}
\end{align*}
and it follows from \eqref{e:Cbeta} that
\begin{align} \label{e:u/V}
\left[ \frac{u}{V(d_D)} \right]_{C^\beta} \leq \frac{c_1 c_8}{r^\beta V(r)} + \frac{c_3 c_4}{V(r)} \leq \frac{c_9}{r^{\beta}V(r)}.
\end{align}

Next, let $x, y \in D$ and let us show that 
\begin{align*}
\lv \frac{u(x)}{V(d_D(x))} - \frac{u(y)}{V(d_D(y))} \rv \leq C \ve x - y \ve^\alpha
\end{align*}
for some $\alpha > 0$. Without loss of generality, we may assume that $r:= d_D(x) \geq d_D(y)$. Fix any $0 < \beta \leq \alpha_1$ and let $p > 1+ \alpha_2/\beta$. If $\ve x - y \ve \leq r^p / 2$, then we have $\ve x - y \ve \leq r/2$ and $y \in B(x, r/2)$ since $r \leq 1$. Thus, by \eqref{e:u/V} we obtain
\begin{align*}
\lv \frac{u(x)}{V(d_D(x))} - \frac{u(y)}{V(d_D(y))} \rv \leq c_9 \frac{\ve x - y \ve^\beta}{r^\beta V(r)} \leq c_{10} \frac{\ve x - y \ve^{\beta - \beta / p}}{V(\ve x - y \ve^{1/p})} \leq c_{11} \ve x - y \ve^{\beta - (\beta + \alpha_2)/p}.
\end{align*}
On the other hand, if $\ve x - y \ve \geq r^p / 2$, let $x_0, y_0 \in \partial D$ be boundary points satisfying $d_D(x) = \ve x - x_0 \ve$ and $d_D(y) = \ve y - y_0 \ve$. Then by the oscillation lemma \ref{l:osc} we have
\begin{align} \label{e:osc1}
\lv \frac{u(x)}{V(d_D)(x)} - \frac{u(x_0)}{V(d_D)(x_0)} \rv \leq c_{12} V(d_D(x))^\gamma, \quad \lv \frac{u(y)}{V(d_D)(y)} - \frac{u(y_0)}{V(d_D)(y_0)} \rv \leq c_{12} V(d_D(y))^\gamma
\end{align}
and 
\begin{align} \label{e:osc2}
\lv \frac{u(x_0)}{V(d_D)(x_0))} - \frac{u(y_0)}{V(d_D)(y_0)} \rv \leq c_{12} V\left( d_D(x) + \ve x - y \ve + d_D(y) \right)^\gamma.
\end{align}
Using inequalities \eqref{e:osc1} and \eqref{e:osc2} we obtain
\begin{align*}
\lv \frac{u(x)}{V(d_D)(x))} - \frac{u(y)}{V(d_D)(y)} \rv \leq c_{12} \left( 2V(r)^\gamma + V(2r + \ve x - y \ve)^\gamma \right) \leq c_{13} \ve x - y \ve^{\alpha_1 \gamma / p}.
\end{align*}
Therefore, taking $\alpha = \min \lb \beta - (\beta + \alpha_2)/p, \alpha_1 \gamma / p \rb$ gives the result.
\qed

\section*{Acknowledgement}

The research of Minhyun Kim and Jaehun Lee is supported by
 the National Research Foundation of Korea (NRF) grant funded by the Korea government (MSIP) : NRF-2016K2A9A2A13003815. 
The research of Panki Kim and Kiahm Lee is supported by the National Research Foundation of Korea(NRF) grant funded by the Korea government(MSIP) 
(No. NRF-2015R1A4A1041675).

\end{document}